\documentclass[reqno,12pt]{amsart} 
\usepackage{amsmath}
\usepackage{amssymb}
\usepackage{amsthm}
\newcommand\NoBlackBoxes{\global\overfullrule0pt}
\NoBlackBoxes
\theoremstyle{plain} 

\def\4{\kern1pt}

\def\6{\vphantom0}

\def\8{\kern-10pt}
\def\7#1{_{(#1)}}

\makeatletter
\let\serieslogo@\relax
\let\@setcopyright\relax

\def\speciallabelmark#1{\def\@currentlabel{#1}}
\makeatother

\begin{document}

\def\ffrac#1#2{\raise.5pt\hbox{\small$\4\displaystyle\frac{\,#1\,}{\,#2\,}\4$}}
\def\ovln#1{\,{\overline{\!#1}}}
\def\ve{\varepsilon}
\def\kar{\beta_r}

\title[Second Order Concentration]{SECOND ORDER CONCENTRATION ON THE SPHERE \\
}

\author{S. G. Bobkov$^{1}$}
\thanks{1) School of Mathematics, University of Minnesota, USA;
Email: bobkov@math.umn.edu}
\address
{Sergey G. Bobkov \newline
School of Mathematics, University of Minnesota  \newline 
127 Vincent Hall, 206 Church St. S.E., Minneapolis, MN 55455 USA
\smallskip}
\email {bobkov@math.umn.edu} 

\author{G. P. Chistyakov$^{2}$}
\thanks{2) Faculty of Mathematics, University of Bielefeld, Germany;
Email: chistyak@math.uni-bielefeld.de}
\address
{Gennadiy P. Chistyakov\newline
Fakult\"at f\"ur Mathematik, Universit\"at Bielefeld\newline
Postfach 100131, 33501 Bielefeld, Germany}
\email {chistyak@math.uni-bielefeld.de}

\author{F. G\"otze$^{3}$}
\thanks{3) Faculty of Mathematics, University of Bielefeld, Germany;
Email: goetze@math.uni-bielefeld.de}
\address
{Friedrich G\"otze\newline
Fakult\"at f\"ur Mathematik, Universit\"at Bielefeld\newline
Postfach 100131, 33501 Bielefeld, Germany}
\email {goetze@mathematik.uni-bielefeld.de}



\begin{abstract}
Sharpened forms of the concentration of measure phenomenon for classes of
functions on the sphere are developed in terms of Hessians of these 
functions.
\end{abstract}

\maketitle
\markboth{S. G. Bobkov, G. P. Chistyakov and F. G\"otze}{Second order
concentration}




\def\theequation{\thesection.\arabic{equation}}
\def\E{{\bf E}}
\def\R{{\bf R}}
\def\C{{\bf C}}
\def\P{{\bf P}}
\def\H{{\rm H}}
\def\Im{{\rm Im}}
\def\Tr{{\rm Tr}}

\def\k{{\kappa}}
\def\M{{\cal M}}
\def\Var{{\rm Var}}
\def\Ent{{\rm Ent}}
\def\O{{\rm Osc}_\mu}

\def\ep{\varepsilon}
\def\phi{\varphi}
\def\F{{\cal F}}
\def\L{{\cal L}}

\def\be{\begin{equation}}
\def\en{\end{equation}}
\def\bee{\begin{eqnarray*}}
\def\ene{\end{eqnarray*}}


\section{{\bf Introduction}}
\setcounter{equation}{0}

\vskip2mm
\noindent
Let $\sigma_{n-1}$ denote the normalized Lebesgue measure on the unit sphere 
$$
S^{n-1} = \{x \in \R^n: |x| = 1\}, \qquad n\geq 2,
$$ 
in the Euclidean $n$-space which is equipped with the canonical inner 
product $\left<\cdot,\cdot\right>$ and the norm $|\cdot|$. The spherical 
concentration phenomenon asserts in particular that mean zero smooth 
functions $f$ on $S^{n-1}$ are of order at most $\frac{1}{\sqrt{n}}$ 
on a large part of the sphere in the sense of $\sigma_{n-1}$. 
This follows already from the Poincar\'e inequality
\be
\int f^2\,d\sigma_{n-1} \leq 
\frac{1}{n-1}\,\int |\nabla_S f|^2\,d\sigma_{n-1},
\en
where $\nabla_S f$ stands for the spherical gradient of $f$. 
Hence, if the integral on the right-hand side is of order 1, the $L^2$-norm 
of $f$ will be of order at most $\frac{1}{\sqrt{n}}$. Moreover, in case 
$|\nabla_S f| \leq 1$, there is a considerably stronger property 
$$
\int e^{(n-1)\, f^2/c}\,d\sigma_{n-1} \leq 2
$$
involving some absolute constant $c>0$. Using a standard normal random 
variable $Z$, it may be stated informally as stochastic dominance
\be
|f| \preceq c\,\frac{|Z|}{\sqrt{n}},
\en
which means a corresponding inequality for the measures/probabilities of 
the tail sets $|f| \geq r$ and $\frac{c}{\sqrt{n}}\, |Z| \geq r$ for all 
$r > 0$.
This property was first emphasized in the early 70's by V. D. Milman 
in the context of the local theory of Banach spaces and led him to the 
understanding of the concentration of measure phenomenon in a much broader 
sense; cf. V. D. Milman, G. Schechtman [M-S], subsequent works by 
M. Talagrand [T1-2] and M. Ledoux [L1-2] for an account of 
basic ideas and results in this direction up to the end of 90's.

Returning to the sphere, in certain problems one deals however with 
smooth functions that turn out to be of a much smaller order than 
$\frac{1}{\sqrt{n}}$. This cannot be guaranteed just by the Lipschitz 
condition $|\nabla_S f| \leq 1$, even if $f$ is orthogonal to linear 
functions in $L^2(S^{n-1},\sigma_{n-1})$ (which play an extremal role 
in (1.1)). Hence, conditions on higher derivatives of $f$ are required.
The aim of this note is to study corresponding conditions in terms 
of the Hessian of $f''_S$ of $f$ by involving both 
the operator norms $\|f_S''(\theta)\|$ and the Hilbert-Schmidt norms
$\|f_S''(\theta)\|_{\rm HS}$ of the matrices $f_S''(\theta)$ 
($\theta \in S^{n-1}$).

Orthogonality of functions on the unit sphere will be understood as 
orthogonality in the Hilbert space $L^2(S^{n-1},\sigma_{n-1})$. 
Restrictions of affine, linear and quadratic functions on $\R^n$ to the
sphere $S^{n-1}$ will be again called affine, 
linear and quadratic functions respectively on the sphere.

\vskip5mm
{\bf Theorem 1.1.} {\it Assume that $f$ is a $C^2$-smooth function on 
$S^{n-1}$ which is orthogonal to all affine functions. If $\|f''_S\| \leq 1$ 
at all points on the sphere and 
$
\int \|f''_S\|_{{\rm HS}}^2\,d\sigma_{n-1} \leq b^2,
$ 
then 
\be
\int \exp\Big\{\frac{n-1}{2(1+b^2)}\, |f|\Big\}\,d\sigma_{n-1} \leq 2.
\en
}

\vskip2mm
By Chebyshev's inequality, (1.3) provides bounds on tails,
which may be written similarly to (1.2) as
$$
|f| \, \preceq \, c_b\Big(\frac{Z}{\sqrt{n}}\Big)^2,
$$
however -- with the right-hand side behaving like $\frac{1}{n}$ 
with respect to the dimension (provided that $b$ is of order 1). 

We refer to Theorem 1.1 as (a variant of) the second order concentration
on the sphere. It is consistent with a second order Poincar\'e-type 
inequality
$$
\int f^2\,d\sigma_{n-1} \leq \frac{1}{2n(n+2)}\,
\int \|f''_S\|_{{\rm HS}}^2\,d\sigma_{n-1},
$$
valid for all smooth $f$ on $S^{n-1}$ that are orthogonal 
to affine functions (with equality attainable for all quadratic
spherical harmonics). This inequality can be derived using the spectral
decomposition of $f$ in spherical harmonics by means of the identity
\be
\int \|f''_S\|_{{\rm HS}}^2\,d\sigma_{n-1} = 
\int f\, \big(\Delta_S (\Delta_S f) + (n-2)\Delta_S f\big)\,d\sigma_{n-1}.
\en
Here and in the sequel $\Delta_S = {\rm Tr}\, f''_S$ denotes the Laplacian 
operator on $S^{n-1}$ which acts diagonally on all homogeneous spherical
harmonics. Although typically $\Delta_S f$ behaves in a more ``chaotic" 
(oscillatory) way than $f$, the average in (1.4) captures and cancels 
such potentially large oscillations.

The conditions on the spherical second derivative in Theorem 1.1 are 
fulfilled, for example, when $\|f''_S\|_{{\rm HS}} \leq b$ on $S^{n-1}$.
However, in applications, one might prefer to deal with functions on the
sphere induced by smooth functions in $\R^n$ or at least in 
a neighbourhood of the sphere via restriction and using the Euclidean 
derivatives of such functions, rather than intrinsic derivatives on $S^{n-1}$. 
Using this Euclidean setup, we may formulate a related statement as follows. 

In the sequel we denote by $f''(x) = \big(\partial_{ij} f(x)\big)_{i,j=1}^n$
the matrix of partial derivatives of $f$ of second order at the point $x$, 
and by $I_n$ the identity $n \times n$ matrix.

\vskip5mm
{\bf Theorem 1.2.} {\it Let $f$ be defined and $C^2$-smooth in some open
neighbourhood of $S^{n-1}$. Assume that it is orthogonal to all affine 
functions and satisfies $\|f'' - a I_n\| \leq 1$ on $S^{n-1}$ together with
\be
\int \|f'' - a I_n\|_{\rm HS}^2\, d\sigma_{n-1} \leq b^2
\en
for some $a \in \R$ and $b \geq 0$. Then
\be
\int \exp\Big\{\frac{n-1}{2(1+4b^2)}\, |f|\Big\}\, d\sigma_{n-1} \leq 2.
\en
}

\vskip2mm
In Theorems 1.1-1.2 one may also start with an arbitrary $C^2$-smooth function 
$f$, but apply the hypotheses and the conclusions (1.3)/(1.6) to the 
projection $Tf$ of $f$ onto the orthogonal complement of the space 
of all affine functions on the sphere in $L^2(S^{n-1},\sigma_{n-1})$. 
The ``affine" part of $f$ may be described as
$l(\theta) = m + \left<v,\theta\right>$ with
$$
m = \int f(x)\,d\sigma_{n-1}(x), \quad
v = n \int xf(x)\,d\sigma_{n-1}(x),
$$
so $Tf(\theta) = f(\theta) - l(\theta)$. For example, if $f$ is even, i.e.
$f(-\theta) = f(\theta)$ for all $\theta \in S^{n-1}$, then
$Tf = f - m$.

In the setting of Theorem 1.2, the functions $Tf$ and $f$ have
identical Euclidean second derivatives. Hence, if we want to obtain an
inequality similar to (1.6) without the orthogonality assumption (still 
assuming conditions on the Euclidean second derivative), we need to
verify that the affine part $l$ is of order $\frac{1}{n}$.
This may be achieved by estimating the $L^2$-norm of $l$ and using the 
well-known fact that the linear functions on the sphere behave like 
Gaussian random variables. If, for definiteness, $f$ has mean zero, then
$$
\|l\|_{L^2}^2 = \frac{1}{n}\,|v|^2 = nI, \quad {\rm where} \ \ 
I = \int\!\!\! \int
\left<x,y\right> f(x)f(y)\,d\sigma_{n-1}(x)d\sigma_{n-1}(y).
$$
Therefore, a natural requirement would be a bound $I \leq \frac{b_0}{n^3}$ 
with $b_0$ of order 1. This leads to a variant of Theorem 1.2 
which is more flexible in applications.

\vskip5mm
{\bf Theorem 1.3.} {\it Let $f$ be defined and $C^2$-smooth in some open
neighbourhood of $S^{n-1}$. Assume that it has mean zero and 
$$
\int\!\!\! \int\left<x,y\right>\,f(x)f(y)\,d\sigma_{n-1}(x)d\sigma_{n-1}(y)
\leq \frac{b_0}{n^3}, \qquad b_0 \geq 0.
$$
If $\|f'' - a I_n\| \leq 1$ holds on $S^{n-1}$ together with $(1.5)$, then
$$
\int \exp\Big\{\frac{n-1}{4(1+b_0^2+4b^2)}\, |f|\Big\}\, d\sigma_{n-1} \leq 2.
$$
}

\vskip2mm
We believe that the second order concentration on the sphere may indeed
be useful in various applications. One motivating example has been the 
question of optimal rates of approximation in the central limit 
theorem for linear forms $X_\theta = \left<X,\theta\right>$, where 
$X = (X_1,\dots,X_n)$ is a given random vector in $\R^n$ whose components 
are not necessarily independent. If the covariance matrix of $X$ has 
a bounded spectral radius, a celebrated result of Sudakov [S] indicates 
that, for $n$ large, the distributions $F_\theta$ of $X_\theta$ are 
concentrated for most of $\theta$ (in the sense of $\sigma_{n-1}$) around 
a certain typical measure $F$ on the real line, which may or may not be 
Gaussian. Many authors studied various aspects of this interesting
phenomenon, and we omit references. 
Let us mention only that one can study the deviations $F_\theta$ 
from $F$ in terms of the Fourier-Stieltjes transforms
$$
f_t(\theta) = \E\,e^{it\left<\theta,X\right>} = 
\int_{-\infty}^{\infty} e^{it\left<\theta,x\right>}\,dF_\theta(x) \qquad
(t \in \R, \ \theta \in \R^n),
$$
which are naturally defined as smooth functions on the whole space 
$\R^n$. By the direct differentiation in $\theta$,
$$
\left<f_t''(\theta)v,w\right> = -t^2\,
\E \left<v,X\right> \left<w,X\right> e^{it\left<\theta,X\right>}.
$$
Here, condition (1.5) leads to a certain correlation-type condition for 
products $X_j X_k$, such that (1.6) will ensure $\frac{1}{n}$-bounds for 
typical deviations of $F_\theta$ from $F$ (in contrast with 
$\frac{1}{\sqrt{n}}$-bounds in the classical Berry-Esseen theorem).
Such improving effects have recently been shown in the work of
B. Klartag and S. Sodin in case of independent summands ([K-S], cf. also [K]). 
As for the general setting, this concentration problem will be dealt 
with in a separate paper and hence will not be discuss it here further.

The proof of Theorems 1.1-1.2 is based on the application of the
logarithmic Sobolev inequality on the sphere and requires derivation 
of bounds on the integrals
\be
\int |\nabla_S f|^2\,d\sigma_{n-1}, \qquad
\int |\nabla f|^2\,d\sigma_{n-1}
\en
in terms of the second derivatives. Basic tools leading to exponential 
bounds under logarithmic Sobolev inequalities are rather universal and 
can be developed in the setting of abstract metric spaces, cf. Section 2.
Then we turn to the case of the sphere and sharpen the Poincar\'e inequality 
by involving the norm $\|f''_S\|$ (Section 3). Sections 5-6 are devoted 
to the estimation of the integrals (1.7). As a preliminary step, 
the identity (1.4) is derived separately in Section 4.
The proofs of Theorem 1.1 and Theorems 1.2-1.3 are completed in Sections 
5 and 7, respectively. After Section 7 we add an Appendix (Sections 8-14)
providing for the readers convenience more details on the underlying 
computations in spherical calculus.

\vskip5mm
{\bf Acknowledgements.} This research was partially supported by NSF grant DMS-1612961, 
the  Humboldt Foundation and SFB 701 at Bielefeld University.
We would like to thank Michel Ledoux for the differential geometric motivation of 
Proposition 4.1, and Bo'az Klartag for the careful reading of the manuscript and
valuable comments.


\vskip10mm
\section{{\bf Logarithmic Sobolev Inequalities on Metric Spaces}}
\setcounter{equation}{0}

\vskip2mm
\noindent
Assume that a metric space $(M,\rho)$ is equipped with a Borel probability 
measure $\mu$. The triple $(M,\rho,\mu)$ is said to satisfy a logarithmic 
Sobolev inequality with constant $\sigma^2 < \infty$, if 
\be
\Ent_\mu(f^2) \leq 2\sigma^2 \int |\nabla f|^2\,d\mu
\en
for any bounded function $f$ on $M$ with finite Lipschitz 
semi-norm $\|f\|_{\rm Lip}$. The optimal value of $\sigma^2$ is then 
called the logarithmic Sobolev constant. 

Here
$$
\Ent_\mu(u) = \int u\log u\,d\mu - \int u\,d\mu \log \int u\,d\mu \qquad
(u \geq 0)
$$
is the entropy functional defined for non-negative measurable
functions on $M$. As for the modulus of the gradient in (2.1), it may be
understood in the generalized sense as
\be
|\nabla f(x)| = \limsup_{y \rightarrow x} \frac{|f(x)-f(y)|}{\rho(x,y)}
\qquad (x \in M).
\en
This function is always Borel measurable, whenever $f$ is continuous.
In this abstract setting, (2.1) actually extends to the larger 
class of all $f$ that have a finite Lipschitz semi-norm on every 
ball in $M$; such functions will be called locally Lipschitz.

Now, define the function
\be
|\nabla^2 f(x)| = |\nabla \, |\nabla f(x)|\,| = 
\limsup_{y \rightarrow x} \frac{|\,|\nabla f(x)|-|\nabla f(y)|\,|}{\rho(x,y)},
\en
which we call a second order modulus of the gradients of $f$.

The Lipschitz property $\|f\|_{\rm Lip} \leq 1$ implies that
$|\nabla f(x)| \leq 1$ for all $x \in M$. The converse is also true, at 
least when $M$ is a (connected) Riemannian manifold. In this case, the 
assumption $|\nabla^2 f(x)| \leq 1$ for every $x$ in $M$ 
means that the function $|\nabla f|$ is Lipschitz. If $|\nabla f|$ is 
locally Lipschitz, then $f$ is of course locally Lipschitz as well.

The next statement indicates how the definition (2.3) could be used 
in applications.

\vskip5mm
{\bf Proposition 2.1.} {\it Assume that a metric probability space 
$(M,\rho,\mu)$ satisfies a logarithmic Sobolev inequality with constant 
$\sigma^2$. Then, for any locally Lipschitz function $f$ on $M$ with 
$\mu$-mean zero, such that $|\nabla f|$ is locally Lipschitz and
$|\nabla^2 f| \leq 1$ on the support of $\mu$, we have
\be
\int \exp\Big\{\frac{1}{2\sigma^2}\, f\Big\}\,d\mu \leq \exp
\Big\{\frac{1}{2\sigma^2} \int |\nabla f|^2\,d\mu\Big\}.
\en
}

\vskip5mm
{\bf Proof.} The argument is based on two general results that relate 
(2.1) to the exponential integrability of Lipschitz functions.
Namely, for any locally Lipschitz $\mu$-integrable function $u$ on $M$,
\be
\int e^{u - \int u\,d\mu}\,d\mu \leq \int e^{\sigma^2 |\nabla u|^2}\,d\mu.
\en
In addition, if $|\nabla u| \leq 1$ on the support of $\mu$, say $M_1$, 
then for all $0 \leq t < \frac{1}{2\sigma^2}$,
\be
\int e^{tu^2}\,d\mu \leq \exp
\Big\{\frac{t}{1 - 2\sigma^2 t}\, \int u^2\,d\mu\Big\}.
\en
On the basis of (2.1), the inequality (2.5) was derived in [B-G], cf. also
[L1-2]. The second inequality, (2.6), is a classical result of Aida, Masuda 
and Shigekawa [A-M-S]. We refer to [B-G] for a detailed discussion. 

We apply (2.6) with $t = \sigma^2\lambda^2$ to the locally Lipschitz
function $u = |\nabla f|$. Since the condition $|\nabla u| \leq 1$
is assumed to hold on $M_1$, we get that
$$
\int e^{\sigma^2\lambda^2 |\nabla f|^2}\,d\mu \leq \exp
\bigg\{\frac{\sigma^2\lambda^2}{1 - 2\sigma^4\lambda^2}\, 
\int |\nabla f|^2\,d\mu\bigg\}, \qquad \lambda^2 < \frac{1}{2\sigma^4}.
$$
On the other hand, since $f$ is locally Lipschitz and has
$\mu$-mean zero, one may apply (2.5), which gives
$$
\int e^{\lambda f}\,d\mu \leq \int e^{\sigma^2 \lambda^2 |\nabla f|^2}\,d\mu.
$$
Hence, the combination of these two bounds yields
$$
\int e^{\lambda f}\,d\mu \leq \exp
\bigg\{\frac{\sigma^2\lambda^2}{1 - 2\sigma^4\lambda^2}\, 
\int |\nabla f|^2\,d\mu\bigg\}.
$$
Here one may choose $\lambda = \frac{1}{2\sigma^2}$, 
and then we arrive at the required inequality (2.4).
\qed

\vskip5mm
When $M$ is an open region in $\R^n$ (with the Euclidean distance), 
the definition (2.1) leads to the usual notion of a logarithmic Sobolev 
inequality, holding for all locally Lipschitz functions on $M$. 
To avoid possible confusion about being locally Lipschitz, let us
emphasize that, when $f$ is differentiable at a given point $x$, 
(2.2) does coincide with the modulus (the length) of the Euclidean gradient.
The same remark applies to the sphere $M = S^{n-1}$ with the geodesic or 
induced Euclidean distances, in which case (2.2) defines $|\nabla_S f(x)|$,
the length of the spherical gradient of $f$. 

The second order modulus of the gradients may also be related to the usual 
(Euclidean) derivatives. Namely, if $f$ is $C^2$-smooth in the open set 
$M$ in $\R^n$, the function $|\nabla f|$ will be locally Lipschitz, and 
\be
|\nabla^2 f(x)| = |\nabla f(x)|^{-1} |f''(x)\nabla f(x)|, \qquad
x \in M.
\en
Here the ratio should be understood as $\|f''(x)\|$ in case 
$|\nabla f(x)| = 0$. In particular,
\be
|\nabla^2 f(x)| \leq \|f''(x)\|.
\en

\vskip5mm
For example, for the quadratic function 
$f(x) = \frac{1}{2}\,\sum_{i=1}^n \lambda_i x_i^2$, $x = (x_1,\dots,x_n)$, 
$$
|\nabla^2 f(x)| \, = \, 
\frac{\sqrt{\sum_{i=1}^n \lambda_i^4 x_i^2}}{\sqrt{\sum_{i=1}^n 
\lambda_i^2 x_i^2}} \, \leq \, \max_i |\lambda_i|.
$$

The identity (2.7) is easily obtained by the direct differentiation.
Thus, in the Euclidean setup Proposition 2.1 may be simplified by using 
the inequality (2.8) as follows.

\vskip5mm
{\bf Corollary 2.2.} {\it Let a probability measure $\mu$ on $\R^n$
satisfy a logarithmic Sobolev inequality with constant $\sigma^2$,
and let a function $f$ be $C^2$-smooth in 
an open neighbourhood of the support of $\mu$. If it has $\mu$-mean zero 
and $\|f''\| \leq 1$ on the support of $\mu$, then
$$
\int \exp\Big\{\frac{1}{2\sigma^2}\, f\Big\}\,d\mu \leq \exp
\Big\{\frac{1}{2\sigma^2} \int |\nabla f|^2\,d\mu\Big\}.
$$
}


\vskip10mm
\section{{\bf Logarithmic Sobolev Inequality on the Sphere}}
\setcounter{equation}{0}

\vskip2mm
\noindent
An important result due to Mueller and Weissler [M-W] sharpens the 
Poincar\'e inequality (1.1) in terms of the logarithmic Sobolev inequality. 
Namely, the logarithmic Sobolev constant of the unit sphere 
$S^{n-1}$, which is equipped with the geodesic metric $\rho$ and 
the uniform measure $\sigma_{n-1}$, coincides with 
the Poincar\'e constant $\sigma^2 = \frac{1}{n-1}$.
That is, for any $C^1$-smooth function $f:S^{n-1} \rightarrow \R$,
\be
\Ent_{\sigma_{n-1}}(f^2) \leq \frac{2}{n-1} \int |\nabla_S f|^2\,d\sigma_{n-1}.
\en 
To see the connection of (3.1) with the concentration phenomenon on the 
sphere in the form (1.2), one may apply (2.6) with $u=f$ and $t = \frac{n-1}{4}$. 

We are also in the position to apply the abstract Proposition 2.1 to 
$(S^{n-1},\rho,\sigma_{n-1})$ and thus involve the second order modulus 
of the gradients, $|\nabla_S^2 f|$. On the unit sphere it is defined according 
to (2.2)-(2.3) with $|\nabla f|$ replaced by
$$
|\nabla_S f(\theta)| 
 \, = \, \limsup_{\theta' \rightarrow \theta} \frac{|f(\theta)-f(\theta')|}{\rho(\theta,\theta')} \qquad 
(\theta, \theta' \in S^{n-1}).
$$
Note that both the geodesic and Euclidean metrics on $S^{n-1}$ may 
equivalently be used for computing the modulus of the gradient of first
and second orders.

For example, the Euclidean derivatives of the linear function 
$f(x) = \left<v,x\right>$ are just $\nabla f(x) = v$ and $f'' = 0$. 
As for the first and second order modulus of its spherical gradient, 
we have
$$
|\nabla_S f(\theta)| = \sqrt{|v|^2 - \left<v,\theta\right>^2} \qquad
(|v|=1),
$$
and, by the chain rule,
\bee
\nabla_S |\nabla_S f(\theta)| 
 & = &
-\frac{1}{2\sqrt{|v|^2 - \left<v,\theta\right>^2}}\
\nabla_S \left(\left<v,\theta\right>^2\right) \\
 & = &
-\frac{1}{\sqrt{|v|^2 - \left<v,\theta\right>^2}}\
\left<v,\theta\right> \nabla_S \left<v,\theta\right> \qquad
(\theta \neq v).
\ene
Hence, $|\nabla_S^2 f(\theta)| = |\left<v,\theta\right>|$
in contrast with $|\nabla^2 f(\theta)| = 0$.

\vskip2mm
To simplify the condition $|\nabla_S^2 f| \leq 1$, one may use the 
following equality which is a full analog of the formula (2.7) mentioned 
before for the case of open regions in $\R^n$.

\vskip5mm
{\bf Lemma 3.1.} {\it Given a $C^2$-smooth function $f$ on $S^{n-1}$, 
$|\nabla_S f|$ has a finite Lipschitz semi-norm and, for all 
$\theta \in S^{n-1}$,
$$
|\nabla_S^2 f(\theta)| \, = \,
|\nabla_S f(\theta)|^{-1}\, |f_S''(\theta)\nabla_S f(\theta)|,
$$
where the right-hand side is understood as $\|f_S''(\theta)\|$ in case 
$|\nabla_S f(\theta)| = 0$. In particular,
$|\nabla_S^2 f(\theta)| \leq \|f_S''(\theta)\|$.
}

\vskip5mm
The proof is given in Appendix (Section 10).

Thus, in order to bound exponential moments of $f$ similarly to (2.4), 
one may require the condition $\|f_S''\| \leq 1$. There is however 
an alternative way based on the application of Corollary 2.2;
the latter would allow us to work with Euclidean derivatives. 
Let us state both consequences of the logarithmic Sobolev inequality (3.1).
Henceforth we shall always understand the mean of functions on the unit 
sphere to be taken with respect to the measure $\sigma_{n-1}$.

\vskip5mm
{\bf Corollary 3.2.} {\it Let $f$ be a $C^2$-smooth function on $S^{n-1}$
with mean zero. If $\|f_S''\| \leq 1$, then
\be
\log
\int \exp\Big\{\frac{n-1}{2}\, f\Big\}\,d\sigma_{n-1} \leq 
\frac{n-1}{2} \int |\nabla_S f|^2\,d\sigma_{n-1}.
\en
Moreover, if $f$ is $C^2$-smooth in an open neighbourhood of the unit 
sphere with $\|f''\| \leq 1$ on $S^{n-1}$, then
\be
\log
\int \exp\Big\{\frac{n-1}{2}\, f\Big\}\,d\sigma_{n-1} \leq 
\frac{n-1}{2} \int |\nabla f|^2\,d\sigma_{n-1}.
\en
}

Applying (3.2) to functions $\ep f$ with $\ep \rightarrow 0$, this 
inequality returns us to (1.1) with an additional factor 2. The condition 
$\|\ep f_S''\| \leq 1$ is fulfilled for all $\ep$ small enough, so any 
constraint on the second derivative may be removed from the conclusion. 
In this sense, Corollary 3.2 provides a sharper form of the Poincar\'e 
inequality.


\vskip10mm
\section{{\bf Second Derivative and Laplacian}}
\setcounter{equation}{0}

\vskip2mm
\noindent
In order to estimate the integral appearing on the right-hand side in (3.2), 
we first derive the formula (1.4), involving the square of the spherical 
Laplacian, i.e. the operator $\Delta_S^2 f = \Delta_S\, \Delta_S f$. Given 
a point $\theta \in S^{n-1}$, it will be convenient to work with the 
spherical second derivative $f''_S(\theta)$ as a symmetric $n \times n$ 
matrix, i.e. as a linear operator on $\R^n$, rather than as a linear operator 
on the tangent space $\theta^\perp$. More precisely, we extend the usual 
Hessian of $f$ at $\theta$ to the whole space by putting
$f''_S(\theta)\theta = 0$ (in particular, both the operator norm and the 
Hilbert-Schmidt norm will not increase for the extended matrix).
The extended Hessian $f_S''(\theta)$ may also be defined as the $n \times n$ 
matrix $B$ with the smallest Hilbert-Schmidt norm, satisfying the Taylor 
expansion
\bee
f(\theta') 
 & = &
f(\theta) + \left<\nabla_S f(\theta),\theta' - \theta\right> \\
 & & + \
\frac{1}{2}\, \left<B(\theta' - \theta),\theta' - \theta\right> +
o\big(|\theta' - \theta|^2\big) \qquad
(\theta' \rightarrow \theta, \ \ \theta' \in S^{n-1}).
\ene

When $f$ is $C^2$-smooth in an open region containing 
the unit sphere, the spherical second derivative is related to 
the Euclidean derivatives by
$$
f_S''(\theta) = P_{\theta^\perp} B P_{\theta^\perp}, \qquad 
B = f''(\theta) - \left<\nabla f(\theta),\theta\right> I_n,
$$
where $P_{\theta^\perp}$ is the projection operator from $\R^n$ 
to the space $\theta^\perp$ orthogonal to $\theta$. Also, recall that 
$\nabla_S f(\theta) = P_{\theta^\perp} \nabla f(\theta)$.

\vskip5mm
{\bf Proposition 4.1.} {\it For any $C^4$-smooth function $f$ on $S^{n-1}$,
\be
\int \|f_S''\|_{\rm HS}^2\,d\sigma_{n-1} \, = \,
\int f\, \big(\Delta^2_S f + (n-2)\, \Delta_S f\big)\, d\sigma_{n-1}.
\en
}

\vskip2mm
One can give a short proof of (4.1) on the basis of the Bochner-Lichnerowicz 
formula in Riemannian Geometry (cf. Remark 4.6 below). Nevertheless, for 
the reader's convenience, we shall provide a direct argument based on 
integration formulas in the multivariate calculus on the sphere which we 
supply in the Appendix, sections A-G. 
The first of these formulas connects the spherical second derivative with 
the iteration of spherical derivatives. The second one is a formula for the 
commutator of the Laplacian and the gradient.

\vskip5mm
{\bf Lemma 4.2.} {\it Given a $C^2$-smooth function $f$ on $S^{n-1}$,
for all $\theta \in S^{n-1}$ and $v \in \R^n$,
\be
f''_S(\theta) v \, = \, \nabla_S \left<\nabla_S f(\theta),v\right> + 
\left<v,\theta\right> \nabla_S f(\theta).
\en
}

\vskip2mm
{\bf Lemma 4.3.} {\it Given a $C^3$-smooth function $f$ on $S^{n-1}$,
for all $\theta \in S^{n-1}$ and $v \in \R^n$, 
$$
\Delta_S \left<\nabla_S f(\theta),v\right> - 
\left<\nabla_S \Delta_S f(\theta),v\right> \, = \,  
(n-3)\left<\nabla_S f(\theta),v\right> - 2\left<v,\theta\right>
\Delta_S f(\theta).
$$
}

The spherical Laplacian appears, in particular, in the integral formula
\be
\int \left<\nabla_S f,\nabla_S g\right>\,d\sigma_{n-1} = -
\int f \Delta_S g\,d\sigma_{n-1}.
\en
The following analogous identity involves a linear weight (cf. Proposition 14.2).

\vskip5mm
{\bf Lemma 4.4.} {\it For all $C^2$-smooth functions $f, g$ on 
$S^{n-1}$ and for any $v \in \R^n$, 
\bee
\int \left<\nabla_S f(\theta),\nabla_S g(\theta)\right>\left<v,\theta\right> 
d\sigma_{n-1}(\theta)
 & = &
-
\int f(\theta)\Delta_S g(\theta)\left<v,\theta\right> d\sigma_{n-1}(\theta) \\
 & & - \ 
\int f(\theta)\left<\nabla_S g(\theta),v\right>\,d\sigma_{n-1}(\theta).
\ene
}

Finally, let us mention how to relate the spherical Laplacian to the Euclidean 
derivatives. The next representation is derived in Section 11, cf. Lemma 11.2;
it will be used in Section 6 in the proof of Theorem 1.2.

\vskip5mm
{\bf Lemma 4.5.} {\it If $f$ is $C^2$-smooth in an open region containing 
the unit sphere, then for any $\theta \in S^{n-1}$, 
$$
\Delta_S f(\theta) = \Delta f(\theta) - (n-1) \left<\nabla f(\theta),\theta\right> 
- \left<f''(\theta)\theta,\theta\right>.
$$
}

\vskip2mm
{\bf Proof of Proposition 4.1.} Using (4.2), one may write
\bee
\int \|f_S''\|_{\rm HS}^2\,d\sigma_{n-1}
 & = &
n \int\!\!\!\int |f_S''(\theta)v|^2\,d\sigma_{n-1}(\theta)d\sigma_{n-1}(v) \\
 & = &
n \int\!\!\!\int \big|\nabla_S \left<\nabla_S f(\theta),v\right> + 
\left<v,\theta\right> \nabla_S f(\theta)\big|^2\,
d\sigma_{n-1}(\theta)d\sigma_{n-1}(v) \\
 & = &
n\,(I_1 + 2 I_2 + I_3),
\ene
where
\bee
I_1
 & = &
\int\!\!\!\int |\nabla_S \left<\nabla_S f(\theta),v\right>|^2\,
d\sigma_{n-1}(\theta)d\sigma_{n-1}(v), \\ 
I_2
 & = &
\int\!\!\!\int \left<\nabla_S \left<\nabla_S f(\theta),v\right>,
\nabla_S f(\theta)\right>
\left<v,\theta\right> d\sigma_{n-1}(\theta)d\sigma_{n-1}(v), \\ 
I_3
 & = &
\int\!\!\!\int |\nabla_S f(\theta)|^2 \left<v,\theta\right>^2
d\sigma_{n-1}(\theta)d\sigma_{n-1}(v). \\ 
\ene
Integration over $v$ immediately gives
$$
I_3 = \frac{1}{n} \int |\nabla_S f|^2\,d\sigma_{n-1},
$$
and according to (4.3),
$$
I_1 = -\int\!\!\!\int \varphi_v(\theta) \Delta_S \varphi_v(\theta)\,
d\sigma_{n-1}(\theta)d\sigma_{n-1}(v), \quad {\rm where} \ \ \ 
\varphi_v(\theta) = \left<\nabla_S f(\theta),v\right>.
$$
To continue, we apply Lemma 4.3, so as to develop 
$\Delta_S \varphi_v(\theta)$ and represent the above integral in the form
$$
I_1 = -\big(I_{11} + (n-3)\,  I_{12} -  2I_{13}\big) 
$$
with
\bee
I_{11}
 & = &
\int\!\!\!\int \left<\nabla_S f(\theta),v\right>
\left<\nabla_S \Delta_S f(\theta),v\right>
d\sigma_{n-1}(\theta)d\sigma_{n-1}(v), \\
I_{12}
 & = &
\int\!\!\!\int \left<\nabla_S f(\theta),v\right>^2
d\sigma_{n-1}(\theta)d\sigma_{n-1}(v), \\
I_{13}
 & = &
\int\!\!\!\int \left<\nabla_S f(\theta),v\right>\left<v,\theta\right>
\Delta_S f(\theta)\, d\sigma_{n-1}(\theta)d\sigma_{n-1}(v). \\
\ene
Let us now integrate over $v$ and apply (4.3) with $g = \Delta_S f$ 
to simplify the first equality as
\bee
I_{11} 
 & = &
\frac{1}{n} \int \left<\nabla_S f(\theta),\nabla_S \Delta_S f(\theta)\right>
d\sigma_{n-1}(\theta) \\
 & = &
-\frac{1}{n} \int f\, \Delta_S (\Delta_S f)\, d\sigma_{n-1} \ = \ 
-\frac{1}{n} \int f\, \Delta^2_S f\, d\sigma_{n-1}.
\ene
We also have
$$
I_{12} = \frac{1}{n} \int |\nabla_S f|^2\,d\sigma_{n-1}, \qquad
I_{13} = \frac{1}{n} \int \left<\nabla_S f(\theta),\theta\right>
\Delta_S f(\theta)\,d\sigma_{n-1}(\theta) = 0.
$$
This finally gives
$$
I_1 \, = \, \frac{1}{n} \int f\, \Delta^2_S f\, d\sigma_{n-1} - 
\frac{n-3}{n} \int |\nabla_S f|^2\,d\sigma_{n-1}.
$$

In order to evaluate the integral $I_2$, we apply Lemma 4.4 with the function
$\left<\nabla_S f(\theta),v\right>$ in place of $f$ and with $f$ in place
of $g$. After integration over $\theta$, we obtain the integral over
the remaining variable $v$, namely,
$$
I_2(v) = -
\int \left<\nabla_S f(\theta),v\right>\Delta_S f(\theta)\left<v,\theta\right> 
d\sigma_{n-1}(\theta) -
\int \left<\nabla_S f(\theta),v\right>^2\,d\sigma_{n-1}(\theta).
$$
The subsequent integration over $v$ cancels the first integral, since 
its integrand will contain the inner product
$\left<\nabla_S f(\theta),\theta\right> = 0$ as a factor. As a result,
\bee
I_2 
 & = &
\int I_2(v)\,d\sigma_{n-1}(v) \\ 
 & = &
-\int\!\!\! \int \left<\nabla_S f(\theta),v\right>^2\,d\sigma_{n-1}(\theta)
d\sigma_{n-1}(v) \ = \ 
- \frac{1}{n} \int |\nabla_S f(\theta)|^2\,d\sigma_{n-1}(\theta).
\ene
It remains to collect these formulas and conclude that
$$
n\,(I_1 + 2 I_2 + I_3) \, = \, \int f\, \Delta^2_S f\, d\sigma_{n-1} - 
(n-2) \int |\nabla_S f|^2\,d\sigma_{n-1}.
$$
Here the last integral can also be written as
$-\int f \Delta_S f\,d\sigma_{n-1}$, cf. (4.3).
\qed

\vskip5mm
{\bf Remark 4.6.}
According to the Bochner-Lichnerowicz formula (cf. e.g. [B-G-L], p. 509), 
for any smooth function $f$ on the Riemannian manifold $(M,g)$,
\be
\frac{1}{2}\,\Delta_g(|\nabla f|^2) = 
\langle\nabla f,\nabla (\Delta_g f)\rangle+|\nabla\nabla f|^2 + 
Ric_g(\nabla f,\nabla f),
\en
where $Ric_g(\nabla f,\nabla f)$ is the Ricci curvature of $(M,g)$ evaluated 
at $\nabla f$. The unit sphere $M = S^{n-1}$ in $\R^n$ has a constant 
curvature, namely, in this case
$$
Ric_g(\nabla f,\nabla f)=(n-2)|\nabla_S f|^2.
$$
Hence, integrating (4.4) over the sphere, we get
\be
\int\Big(\frac{1}{2}\, \Delta_S(|\nabla_S f|^2) - 
\langle\nabla_S f,\nabla_S (\Delta_S f)\rangle\Big)\,d\sigma_{n-1} =
\int (|\nabla_S\nabla_S f|^2+(n-2)|\nabla_S f|^2)\,d\sigma_{n-1}.
\en
On the other hand, $\int\Delta_S(|\nabla_S f|^2)\,d\sigma_{n-1} = 0$, 
$$
\int|\nabla_S f|^2\,d\sigma_{n-1}=-\int f\Delta_S f\,d\sigma_{n-1},
$$
(recall (4.3)), and 
$$
\int\langle\nabla_S f,\nabla_S (\Delta_S f)\rangle\,d\sigma_{n-1} =
-\int f \Delta_S^2 f\,d\sigma_{n-1}.
$$
Applying these relations in (4.5), we arrive at (4.1).


\vskip10mm
\section{{\bf Expansions in Spherical Harmonics}}
\setcounter{equation}{0}

\vskip2mm
\noindent
Using Proposition 4.1, one may study relations of the form
\be
c \int |\nabla_S f|^2\,d\sigma_{n-1} \leq
\int \|f''_S\|_{{\rm HS}}^2\,d\sigma_{n-1} \qquad (c>0)
\en
by means of the orthogonal expansion in spherical harmonics,
\be
f = \sum_{d=0}^\infty f_d \qquad (f_d \in H_d).
\en
As is well-known (cf. e.g. [S-W]), the Hilbert space $L^2(S^{n-1})$ 
can be decomposed into a sum of orthogonal linear subspaces $H_d$, 
$d = 0,1,2,\dots$, consisting of all $d$-homogeneous harmonic polynomials
(more precisely - restrictions of such polynomials to the sphere).
Any element $f_d$ of $H_d$ represents an eigenfunction of the Laplacian, 
with the eigenvalue $- d(n+d-2)$. That is,
$$
\Delta_S f_d = - d(n+d-2)\,f_d,
$$
and hence
$$
\Delta^2_S f_d \, = \, d^2(n+d-2)^2\, f_d.
$$
As a result,
$$
\Delta_S f = -\sum_{d=1}^\infty d(n+d-2) f_d, \qquad
\Delta^2_S f = \sum_{d=1}^\infty d^2(n+d-2)^2 f_d
$$
which should be understood as equalities in $L^2$ (Note that both 
$\Delta_S f$ and $\Delta^2_S f$ are continuous functions, as long as
$f$ is $C^4$-smooth).

According to the representation (4.1), (5.1) is equivalent to
\be
\int f\,\Delta^2_S f\, d\sigma_{n-1} \, \geq \, 
-(c + n-2)\int f\,\Delta_S f\, d\sigma_{n-1}.
\en
Moreover, since the spherical harmonics serve as eigenfunctions both for
$\Delta_S$ and $\Delta^2_S$, the last inequality need to be verified
for elements $f_d$ of $H_d$ only. Here, both integrals are vanishing for 
constant functions, i.e. for $f \in H_d$ with $d=0$. If $d \geq 1$, (5.3) becomes
\be
c \, \leq \, d^2 + (d-1)(n-2).
\en
Thus, if we want to involve in (5.1) all $C^2$-smooth functions $f$, the optimal 
value of $c$ is described as the minimum of the right-hand side of (5.4) over all 
$d \geq 1$. The minimum is achieved for $d=1$ which leads to the optimal 
value $c=1$. However, if we require that $f$ is orthogonal to all linear functions, 
it means that we only allow the values $d \geq 2$ in (5.4), and then 
the optimal value is $c = n+2$. As a result, we have proved:

\vskip5mm
{\bf Proposition 5.1.} {\it For any $C^2$-function $f$ on $S^{n-1}$,
$$
\int |\nabla_S f|^2\,d\sigma_{n-1} \leq
\int \|f''_S\|_{{\rm HS}}^2\,d\sigma_{n-1},
$$
where equality is attained for all linear functions. Moreover, if $f$ 
is orthogonal to all linear functions with respect to $\sigma_{n-1}$, 
then 
\be
\int |\nabla_S f|^2\,d\sigma_{n-1} \leq \frac{1}{n+2}\,
\int \|f''_S\|_{{\rm HS}}^2\,d\sigma_{n-1}
\en
with equality attainable for all quadratic harmonics.
}

\vskip5mm
The expansion (5.2) is commonly used to derive Poincar\'e-type inequalities
such as (1.1). If we require additionally that $f$ should be orthogonal 
to all linear functions, the constant will slightly improve only, since then
$$
\int f^2\,d\sigma_{n-1} \leq \frac{1}{2n} \int |\nabla_S f|^2\,d\sigma_{n-1}.
$$
This bound may be combined with (5.5) to get a second order Poincar\'e-type
inequality which was mentioned in the Introduction. But, one can also apply 
(5.2) directly in the representation (4.1). Indeed, on 
spherical harmonics $f_d$ of $H_d$, the inequality of the form
$
c\int f^2\,d\sigma_{n-1} \leq \int \|f''_S\|_{{\rm HS}}^2\,d\sigma_{n-1}
$
becomes
$$
c \, \leq \, d(n+d-2)\,\big(d(n+d-2) - (n-2)\big).
$$
Since the right-hand side is an increasing function of $d$, we arrive at:

\vskip5mm
{\bf Proposition 5.2.} {\it For any $C^2$-function $f$ on $S^{n-1}$ with
mean zero,
\be
\int f^2\,d\sigma_{n-1} \leq \frac{1}{n-1}\,
\int \|f''_S\|_{{\rm HS}}^2\,d\sigma_{n-1},
\en
where equality is attained for all linear functions. Moreover, if $f$ 
is orthogonal to all linear functions with respect to $\sigma_{n-1}$, 
then 
\be
\int f^2\,d\sigma_{n-1} \leq \frac{1}{2n(n+2)}\,
\int \|f''_S\|_{{\rm HS}}^2\,d\sigma_{n-1}
\en
with equality attainable for all quadratic harmonics.
}

\vskip5mm
An interesting consequence of (5.6) is the statement that the equality 
$f''_S = 0$ is possible for constant functions, only (in contrast with the
Euclidean Hessian).

\vskip5mm
{\bf Remark 5.3.} It is much easier to derive (5.7) with suboptimal,
although asymptotically correct constants as $n$ tends to infinity,
without appealing to Proposition 4.1. The argument is based on the double
application of the Poincar\'e inequality (1.1). Orthogonality of $f$ to
all linear functions ensures that the function 
$\theta \rightarrow \left<\nabla_S f(\theta),v\right>$ has mean
zero for any $v \in \R^n$. So, using the identity (4.2), we get
\bee
(n-1)\int \left<\nabla_S f(\theta),v\right>^2\,d\sigma_{n-1}(\theta) 
 & \leq &
\int 
|f''_S(\theta) v - \left<v,\theta\right> 
\nabla_S f(\theta)|^2\,d\sigma_{n-1}(\theta) \\
 & \hskip-20mm = &
\hskip-10mm \int |f''_S(\theta) v|^2\,d\sigma_{n-1}(\theta) + 
\int \left<v,\theta\right>^2 |\nabla_S f(\theta)|^2\,d\sigma_{n-1}(\theta) \\
 & & \hskip-15mm - \
2 \int \left<f_S''(\theta) \nabla_S f(\theta),v\right>\left<v,\theta\right>
\,d\sigma_{n-1}(\theta).
\ene
The next integration over $d\sigma_{n-1}(v)$ cancels the last integral
(due to $f_S''(\theta) \theta = 0$), and we are led to
$$
(n-2)\int |\nabla_S f(\theta)|^2\,d\sigma_{n-1}(\theta) 
 \, \leq \, \int \|f''_S(\theta)\|_{\rm HS}^2\,d\sigma_{n-1}(\theta).
$$
If $f$ has mean zero, the left integral may be estimated from below
according to (1.1), which thus gives
$$
\int f^2\,d\sigma_{n-1} \leq \frac{1}{(n-1)(n-2)}\,
\int \|f''_S\|_{{\rm HS}}^2\,d\sigma_{n-1}, \qquad n \geq 3.
$$
The constant in this inequality is slightly worse than (5.7), and we loose
information about extremal functions.

The above argument is also applicable in the Euclidean setup when dealing 
with a probability measure $\mu$ on $\R^n$ satisfying a Poincar\'e-type
inequality
$$
\int f^2\,d\mu \leq \sigma^2 \int |\nabla f|^2\,d\mu \qquad
\bigg(\int f\,d\mu = 0\bigg).
$$
For example, the standard Gaussian measure with density
$
\frac{d\mu(x)}{dx} = (2\pi)^{-n/2}\,e^{-|x|^2/2}
$
has the Poincar\'e constant $\sigma^2 = 1$, which yields a second order
Poincar\'e-type inequality
$$
\int f^2\,d\mu \leq \frac{1}{2}\,
\int \|f''_S\|_{{\rm HS}}^2\,d\mu.
$$
It holds true in the class of all $C^2$-smooth functions $f$ 
on $\R^n$ that are orthogonal to all affine functions in $L^2(\mu)$.
However, in the general case, orthogonality to linear functions
should be replaced with the requirement $\int \nabla f\,d\mu = 0$.

\vskip5mm
We are now prepared to complete the proof of Theorem 1.1.

\vskip5mm
{\bf Proof of Theorem 1.1.}
Let us return to the bound (3.2) of Corollary 3.2. Using (5.5), we then
get
$$
\log
\int \exp\Big\{\frac{n-1}{2}\, f\Big\}\,d\sigma_{n-1} \, \leq \, 
\frac{n-1}{2(n+2)} \int  \|f''_S\|_{{\rm HS}}^2\,d\sigma_{n-1} \, \leq \,
\frac{1}{2}\,b^2,
$$
and using a similar inequality for the function $-f$,
$$
\int e^{\frac{n-1}{2}\, |f|}\,d\sigma_{n-1} \, \leq \,
\int e^{\frac{n-1}{2}\, f}\,d\sigma_{n-1} +
\int e^{-\frac{n-1}{2}\,f}\,d\sigma_{n-1} \, \leq \, 2e^{b^2/2}.
$$
It follows that, for any $\lambda \geq 1$,
$$
\int e^{\frac{n-1}{2}\, |f|/\lambda}\,d\sigma_{n-1}
 \, \leq \,
\Big(\int e^{\frac{n-1}{2}\, |f|}\,d\sigma_{n-1}\Big)^{1/\lambda} \, \leq \, 
 (2e^{b^2/2})^{1/\lambda}.
$$
It remains to note that $(2e^{b^2/2})^{1/\lambda} = 2$ for 
$\lambda = 1 + \frac{b^2}{\log 4} \leq 1+b^2$.
\qed


\vskip10mm
\section{{\bf Bounds on the $L^2$-Norm of the Euclidean Gradient}}
\setcounter{equation}{0}

\vskip2mm
\noindent
We now turn back to Theorem 1.2 while invoking the second 
bound of Corollary 3.2. Hence, we need an analog of (5.5) for the modulus
of the Euclidean gradient. Assume that a function $f$ is defined and 
$C^2$-smooth in some neighbourhood $G$ of $S^{n-1}$.

\vskip5mm
{\bf Proposition 6.1.} {\it If f is orthogonal to all 
linear functions with respect to $\sigma_{n-1}$, then
\be
\int |\nabla f|^2\,d\sigma_{n-1} \leq \frac{5}{n-1}
\int \|f''\|_{\rm HS}^2\,d\sigma_{n-1}.
\en
}

\vskip2mm
At the of the proof it will be apparent that for growing dimensions
the constant 5 may be asymptotically improved to 2.

\vskip5mm
{\bf Proof.} Since the spherical gradient $\nabla_S f(\theta)$ represents
the projection of the usual gradient $\nabla f(\theta)$ to the subspace
$\theta^\perp$ of $\R^n$ orthogonal to $\theta$, we have
$$
|\nabla f|^2 \, = \, |\nabla_S f(\theta)|^2 + 
\left<\nabla f(\theta),\theta\right>^2.
$$

As a preliminary step, first we show that
\be
\int |\nabla_S f|^2\,d\sigma_{n-1} \leq \frac{1}{n-1}
\int \|f''\|_{\rm HS}^2\,d\sigma_{n-1}.
\en
Write
\be
\int |\nabla_S f|^2\,d\sigma_{n-1} \, = \,
\int |\nabla f(\theta)|^2\,d\sigma_{n-1}(\theta) - 
\int\left<\nabla f(\theta),\theta\right>^2\,d\sigma_{n-1}(\theta)
\en
and represent
\be
\int |\nabla f|^2\,d\sigma_{n-1} = n
\int\!\!\!\int \left<\nabla f(\theta),v\right>^2\,d\sigma_{n-1}(\theta)
d\sigma_{n-1}(v).
\en

The assumption that $f$ is orthogonal to all linear functions is equivalent
to the property that every function of the form
$$
\left<\nabla_S f(\theta),v\right>  =  
\left<\nabla f(\theta),v\right> - 
\left<\nabla f(\theta),\theta\right> \left<v,\theta\right>
$$ 
has $\sigma_{n-1}$-mean zero (cf. Proposition 14.1). Hence
$$
\int \left<\nabla f(\theta),v\right>\,d\sigma_{n-1}(\theta) =
\int 
\left<\nabla f(\theta),\theta\right> \left<v,\theta\right>\,d\sigma_{n-1}(\theta),
$$
and, by the Cauchy-Schwarz inequality,
\be
\Big(\int \left<\nabla f(\theta),v\right>\,d\sigma_{n-1}(\theta)\Big)^2
\leq \frac{1}{n} \int 
\left<\nabla f(\theta),\theta\right>^2\,d\sigma_{n-1}(\theta).
\en

To estimate the $L^2$-norm of $\left<\nabla f(\theta),v\right>$, 
one may apply the Poincar\'e inequality (1.1). Since
$u(x) = \left<\nabla f(x),v\right>$ has gradient 
$\nabla u(x) = f''(x)v$, we have, by (6.5),
$$
\int \left<\nabla f(\theta),v\right>^2\,d\sigma_{n-1}(\theta) \leq 
\frac{1}{n} \int 
\left<\nabla f(\theta),\theta\right>^2\,d\sigma_{n-1}(\theta) +
\frac{1}{n-1}\,\int |f''(\theta)v|^2\,d\sigma_{n-1}(\theta).
$$
Using this bound in (6.4) and integrating over $v$, we get
$$
\int |\nabla f|^2\,d\sigma_{n-1} \, \leq \,
\int 
\left<\nabla f(\theta),\theta\right>^2\,d\sigma_{n-1}(\theta) +
\frac{1}{n-1}\,\int \|f''(\theta)\|^2_{\rm HS}\,d\sigma_{n-1}(\theta).
$$
It remains to insert this bound in (6.3) which gives (6.2).

Now, rewrite (6.3) as
\be
\int |\nabla f|^2\,d\sigma_{n-1} =  
\int |\nabla_S f|^2\,d\sigma_{n-1} +
\int \left<\nabla f(\theta),\theta\right>^2\,d\sigma_{n-1}(\theta).
\en
Here, the first integral on the right-hand side is estimated in terms 
of $\|f''\|_{\rm HS}^2$ by (6.2), and our next task will be to derive 
a suitable bound on the $L^2$-norm of the function 
$\left<\nabla f(\theta),\theta\right>$. 
To this aim, we employ the representation of Lemma 4.5 for
the spherical Laplacian in terms of the Euclidean derivatives.
Since in general (by (4.3)),
$$
\int \Delta_S f\,d\sigma_{n-1} = -
\int \left<\nabla_S 1,\nabla_S f\right> d\sigma_{n-1} = 0,
$$ 
Lemma 4.5 yields
\be
(n-1) \int \left<\nabla f(\theta),\theta\right> d\sigma_{n-1}(\theta) =
\int \big(\Delta f(\theta) - 
\left<f''(\theta)\theta,\theta\right>\big)\, d\sigma_{n-1}(\theta).
\en
Here the second integrand is equal to
$$
I = \sum_{i,j = 1}^n \partial_{ij} f(\theta) a_{ij} \quad {\rm with} \ \
a_{ij} = \delta_{ij} - \theta_i \theta_j.
$$
Note that
$$
\sum_{i,j=1}^n a_{ij}^2 \, = \,
\sum_{i \neq j}^n \theta_i^2 \theta_j^2 + \sum_{i=1}^n (1 - \theta_i^2)^2 \, = \,
1 + \sum_{i=1}^n \big((1 - \theta_i^2)^2 - \theta_i^4\big) \, = \, n-1.
$$
Hence, by Cauchy's inequality,
$$
I^2 \leq \sum_{i,j=1}^n (\partial_{ij} f(\theta))^2 \sum_{i,j=1}^n a_{ij}^2 =
(n-1)\,\|f''(\theta)\|_{{\rm HS}}^2,
$$
and by another application of the Cauchy-Schwarz inequality in (6.7),
\be
\bigg(\int \left<\nabla f(\theta),\theta\right> d\sigma_{n-1}(\theta)\bigg)^2
 \, \leq \, \frac{1}{n-1} 
 \int \|f''\|_{\rm HS}^2\,d\sigma_{n-1}.
\en

Next, consider the function $u(x) = \left<\nabla f(x),x\right>$ and
restrict its gradient $\nabla u(x) = \nabla f(x) + f''(x)x$ to the unit 
sphere. Projecting it to $\theta^\perp$, we obtain the spherical gradient
$$
\nabla_S u(\theta) = \nabla_S f(\theta) + 
P_{\theta^\perp} \big(f''(\theta)\theta\big), \qquad \theta \in S^{n-1}.
$$
In particular, by the triangle inequality,
$$
|\nabla_S u(\theta)| \leq |\nabla_S f(\theta)| + \|f''(\theta)\|.
$$
Furthermore, the square of the right-hand side can be estimated by using 
the elementary inequality 
$(x+y)^2 \leq \frac{\lambda}{\lambda-1}\, x^2 + \lambda y^2$ 
($x,y \geq 0$, $\lambda>1$), which implies
$$
|\nabla_S u(\theta)| \leq
\frac{\lambda}{\lambda-1}\,|\nabla_S f(\theta)|^2 + \lambda\,\|f''(\theta)\|^2.
$$
Hence, using the Poincar\'e inequality together with (6.8), 
and increasing the operator norm to the Hilbert-Schmidt norm, we get
\bee
\int u^2\,d\sigma_{n-1} 
 & \leq &
\Big(\int u\, d\sigma_{n-1}\Big)^2 + \frac{1}{n-1}
\int |\nabla_S u|^2\, d\sigma_{n-1} \\
 & \leq &
\frac{1}{n-1} \int \|f''\|_{\rm HS}^2\,d\sigma_{n-1}+
\frac{1}{n-1} \int \left(\frac{\lambda}{\lambda-1}\,
|\nabla_S f|^2 + \lambda\,\|f''\|_{\rm HS}^2\right)\, d\sigma_{n-1}.
\ene
Thus,
$$
(n-1) \int \left<\nabla f(\theta),\theta\right>^2\,d\sigma_{n-1}(\theta) 
 \, \leq \,
\frac{\lambda}{\lambda-1} \int |\nabla_S f|^2\, d\sigma_{n-1}
+ (\lambda+1) \int \|f''\|_{\rm HS}^2\,d\sigma_{n-1}.
$$

It remains to return to (6.6) and combine the above bound with (6.2).
Adding and collecting the coefficients, it gives
$$
(n-1) \int |\nabla f|^2\, d\sigma_{n-1} \, \leq \,
\Big(\frac{1}{n-1} \ \frac{\lambda}{\lambda-1} + \lambda + 1\Big) 
\int \|f''\|_{\rm HS}^2\,d\sigma_{n-1}.
$$
The quantity $\frac{1}{n-1}\,\frac{\lambda}{\lambda-1} + \lambda + 1$ 
is minimized at $\lambda = 1 + \frac{1}{\sqrt{n-1}}$, which leads to
\be
\int |\nabla f|^2\,d\sigma_{n-1} \leq \frac{c_n}{n-1}
\int \|f''\|_{\rm HS}^2\,d\sigma_{n-1}, \qquad
c_n = 1 + \Big(1 + \frac{1}{\sqrt{n-1}}\Big)^2.
\en
Clearly, $c_n \leq 5$, thus proving (6.1).
\qed

\vskip5mm
Note that $c_n \rightarrow 2$ as $n \rightarrow \infty$.
So, the constant 5 in (6.1) may be improved for large values of $n$.

Combining (6.1) with the Poincar\'e inequality (1.1), we get
a second order Poincar\'e-type inequality in the Euclidean setup,
$$
\int (f-m)^2\,d\sigma_{n-1} \leq \frac{5}{(n-1)^2}
\int \|f''\|_{\rm HS}^2\,d\sigma_{n-1},
$$
assuming that $f$ is orthogonal to all linear functions, and where $m$ 
is the mean of $f$ with respect to $\sigma_{n-1}$. Here the left integral 
will not change when it is applied to 
$f_a(x) = f(x) - \frac{a}{2}\,|x|^2$ in place of $f$, 
while the right integral will depend on $a$. More precisely, we get
$$
\int (f-m)^2\,d\sigma_{n-1} \leq \frac{5}{(n-1)^2}
\int \|f'' - aI_n\|_{\rm HS}^2\,d\sigma_{n-1}.
$$
Hence, we arrive at:

\vskip5mm
{\bf Corollary 6.2.} {\it If $f$ is orthogonal to all 
affine functions with respect to $\sigma_{n-1}$, then for any $a \in \R$,
$$
\int f^2\,d\sigma_{n-1} \leq \frac{5}{(n-1)^2}
\int \|f'' - aI_n\|_{\rm HS}^2\,d\sigma_{n-1}.
$$
}


\vskip10mm
\section{{\bf Proof of Theorems 1.2-1.3}}
\setcounter{equation}{0}

\vskip2mm
\noindent
Having proved Proposition 6.1, the proof of Theorem 1.2 is almost identical 
to the proof of Theorem 1.1. 

\vskip5mm
{\bf Proof of Theorem 1.2.}
Let $f$ be orthogonal to all affine functions with mean $m$. Applying (6.1) 
to the function $f-m$ in the bound (3.3) of Corollary 3.2, we get
$$
\log
\int \exp\Big\{\frac{n-1}{2}\, (f-m)\Big\}\,d\sigma_{n-1} \, \leq \, 
\frac{5}{2} \int  \|f''\|_{{\rm HS}}^2\,d\sigma_{n-1}.
$$
Applying it to $f_a(x) = f(x) - \frac{a}{2}\,|x|^2$ in place of $f$, we get
$$
\log
\int \exp\Big\{\frac{n-1}{2}\, (f-m)\Big\}\,d\sigma_{n-1} \, \leq \, 
\frac{5}{2} \int  \|f'' - aI_n\|_{{\rm HS}}^2\,d\sigma_{n-1}
 \, \leq \, \frac{5}{2}\,b^2.
$$
Assuming that $m=0$ and applying a similar inequality to the function $-f$,
we obtain
$$
\int e^{\frac{n-1}{2}\, |f|}\,d\sigma_{n-1} \, \leq \, 2 e^{5 b^2/2}.
$$
Hence, for any $\lambda \geq 1$,
$$
\int e^{\frac{n-1}{2}\, |f|/\lambda}\,d\sigma_{n-1}
 \, \leq \,
\Big(\int e^{\frac{n-1}{2}\, |f|}\,d\sigma_{n-1}\Big)^{1/\lambda} \, \leq \, 
 (2e^{5b^2/2})^{1/\lambda}.
$$
It remains to note that $(2e^{5b^2/2})^{1/\lambda} = 2$
for $\lambda = 1 + \frac{5b^2}{\log 4} \leq 1 + 3.7\,b^2$.
\qed

\vskip5mm
{\bf Proof of Theorem 1.3.} Let $l(\theta) = \left<v,\theta\right>$ be
the linear part of $f$, and recall that
$$
|v|^2 = n^2 I, \qquad
I = \int\!\!\! \int
\left<x,y\right> f(x)f(y)\,d\sigma_{n-1}(x)d\sigma_{n-1}(y).
$$
To control Gaussian tails of $l$ under $\sigma_{n-1}$, we apply 
an exponential bound
$$
\int e^{t l(\theta)}\,d\sigma_{n-1}(\theta) \, \leq \,
e^{\frac{t^2}{2(n-1)}\,|v|^2}, \qquad t \in \R,
$$
which is implied by the logarithmic Sobolev inequality on the sphere, (3.1).
Choosing $t = n-1$ and using the assumption $I \leq \frac{b_0}{n^3}$, we get
$\int e^{(n-1) |l|}\,d\sigma_{n-1} \leq 2 e^{b_0^2/2}$ and hence
$$
\int \exp\Big\{\frac{n-1}{1+b_0^2}\, |l|\Big\}\, d\sigma_{n-1} \leq 2.
$$

On the other hand, by Theorem 1.2 with the same assumption on the second
derivative of $f$, we have
$$
\int \exp\Big\{\frac{n-1}{2(1+4b^2)}\, |Tf|\Big\}\, d\sigma_{n-1} \leq 2.
$$
Using $|f| \leq |Tf| + |l|$ and applying the Cauchy-Schwarz inequality,
we conclude that
$$
\int e^{(n-1)\, |f|/2\lambda}\,d\sigma_{n-1}
 \, \leq \,
\Big(\int e^{(n-1)\, |Tf|/\lambda}\,d\sigma_{n-1}\Big)^{1/2}
\Big(\int e^{(n-1)\, |l|/\lambda}\,d\sigma_{n-1}\Big)^{1/2} \, \leq \, 2,
$$
provided that $\lambda \geq 2(1+4b^2)$ and $\lambda \geq 1 + b_0^2$.
\qed

\vskip10mm
\section{{\bf Appendix A. Definitions of Spherical Derivatives}}
\setcounter{equation}{0}

\vskip2mm
\noindent
A function $f$ defined on the unit sphere $S^{n-1}$ is $C^p$-smooth, 
$p = 1,2,\dots$, if it can be extended to some open set
containing $S^{n-1}$ as a $C^p$-smooth function (in the usual sense). 
This is one of the well-known definitions of smoothness on the sphere.

If $f$ is $C^1$-smooth on $S^{n-1}$, then at every point $\theta \in S^{n-1}$ 
it admits the Taylor expansion up to the linear term
\be
f(\theta') = f(\theta) + \left<v,\theta' - \theta\right> +
o\big(|\theta' - \theta|\big), \quad {\rm as} \ \ 
\theta' \rightarrow \theta, \ \ \theta' \in S^{n-1},
\en
with some $v \in \R^n$. If $v$ has the smallest length (Euclidean norm) 
among all such vectors, it is called the spherical derivative or 
gradient of $f$ at $\theta$ and is denoted $\nabla_S f(\theta)$. 

This notion of the derivative of $f$ is independent of the choice of a 
smooth extension of $f$ in an open neighbourhood of the sphere in $\R^n$.
If $f$ is $C^1$-smooth in a neighbourhood of the unit 
sphere, then (8.1) holds with the usual (Euclidean) gradient 
$v = \nabla f(\theta)$, and the spherical gradient may be described as 
\bee
\nabla_S f(\theta) 
 & = & 
P_{\theta^\perp} \nabla f(\theta) \\
 & = &
\nabla f(\theta) - \left<\nabla f(\theta),\theta\right> \theta, 
\ene
where $P_{\theta^\perp}$ is the (orthogonal) projection operator 
from $\R^n$ to $\theta^\perp$ (the tangent space). 

In particular, $\left<\nabla_S f(\theta),\theta\right> = 0$ and
$|\nabla_S f(\theta)| \leq |\nabla f(\theta)|$ for any 
$\theta \in S^{n-1}$.

The spherical gradient of any $C^1$-function represents a continuous 
vector-valued function on $S^{n-1}$.

Analogously (as was already stressed in Section 4), the second derivative 
of any $C^2$-smooth function $f$ on the unit sphere at a given point 
$\theta \in S^{n-1}$ may be introduced via a Taylor expansion up to 
the quadratic term
\be
f(\theta')  \, = \, f(\theta) + 
\left<\nabla_S f(\theta),\theta' - \theta\right> +
\frac{1}{2} \left<B(\theta' - \theta),\theta' - \theta\right> +
o\big(|\theta' - \theta|^2\big),
\en
where $\theta' \rightarrow \theta$, $\theta' \in S^{n-1}$, and
$B$ is some $n \times n$ matrix (with real entries).

Recall that the space ${\bf M}_n$ of all $n \times n$ 
matrices is naturally identified with the Euclidean space $\R^{n \times n}$ 
with its inner product and the Euclidean norm
$$
\|B\|_{\rm HS} = \bigg(\sum_{i,j = 1}^n B_{ij}^2\bigg)^{1/2}
$$
called the Hilbert-Schmidt norm of $B$. The collection of all $B$ satisfying 
(8.2) represents an affine subspace of ${\bf M}_n$. Therefore, among all
of them, there exists a unique matrix which has the smallest Hilbert-Schmidt 
norm. It can be called the (spherical) second derivative of $f$ at the 
point $\theta$ and will be denoted $f_S''(\theta)$.

If $f$ is $C^2$-smooth in an open neighborhood of $S^{n-1}$, then in 
accordance with the usual Taylor expansion, (8.2) holds with the matrix
$
B_0 = f''(\theta) - \left<\nabla f(\theta),\theta\right> I_n.
$
More generally, given $A \in {\bf M}_n$, the matrix $B_0 - A$ satisfies 
(8.2), if and only if
$$
\left<A(\theta' - \theta),\theta' - \theta\right> = 
o\big(|\theta' - \theta|^2\big) \qquad
(\theta' \rightarrow \theta, \ \ \theta' \in S^{n-1}).
$$
But this is equivalent to saying that 
$\left<Ax,x\right> = 0$ for all $x \in \theta^\perp$. This condition 
defines a linear subspace $L$ of ${\bf M}_n$, and the problem 
$$
\|B_0 - A\|_{{\rm HS}} \rightarrow \min \ {\rm over \ all} \ A \in L
$$
is then solved uniquely for $B = B_0 - A$ being the orthogonal projection 
in ${\bf M}_n$ of $B_0$ to the linear space $L^\perp$ of all matrices
orthogonal to $L$.
In fact, since $B_0$ is symmetric, in this minimization problem one may 
restrict ourselves to symmetric matrices, and by a simple algebra, we 
arrive at the following description.

\vskip5mm
{\bf Proposition 8.1.} {\it The spherical second derivative of $f$ at each
point $\theta \in S^{n-1}$ is a symmetric matrix, which is given by the 
orthogonal projection
$$
f_S''(\theta) = P_{L_{\theta}^\perp} B, \qquad
B = f''(\theta) - \left<\nabla f(\theta),\theta\right> I_n,
$$
to the orthogonal complement of the linear subspace $L_{\theta}$ of 
all symmetric matrices $A$ in ${\bf M}_n$ such that $Ax = 0$ 
for all $x \in \theta^\perp$. Equivalently,
$$
f_S''(\theta) = P_{\theta^\perp} B P_{\theta^\perp}.
$$
}

One immediate consequence of this description is that 
$f_S''(\theta)\theta = 0$ and hence the vectors $f_S''(\theta)v$ are 
orthogonal to $\theta$, for all $\theta \in S^{n-1}$ and $v \in \R^n$.

One should also emphasize the contraction property
$$
\|f_S''(\theta)\|_{{\rm HS}} \, \leq \, 
\|f''(\theta) - \left<\nabla f(\theta),\theta\right> I_n\|_{{\rm HS}}
$$
and similarly for the operator norm.

\vskip10mm
\section{{\bf Appendix B. Second Order Gradients}}
\setcounter{equation}{0}

\vskip2mm
\noindent
Let us now turn to Lemma 4.2 with its identity
\be
f''_S(\theta) v = \nabla_S \left<\nabla_S f(\theta),v\right> + 
\left<v,\theta\right> \nabla_S f(\theta).
\en
Note that the usual first and second derivatives are connected by
\be
f''(x) v = \nabla \left<\nabla f(x),v\right> \qquad (v \in \R^n).
\en
As follows from (9.1), we have a similar property for the spherical 
derivatives -- however for $v$ in the tangent space, only.

\vskip5mm
{\bf Proof of Lemma 4.2.} We may assume that $f$ is defined and 
$C^2$-smooth on an open subset $G$ of $\R^n$ containing the unit sphere.
To compute the spherical gradient for the function
$$
\psi(\theta) \, = \, \left<\nabla_S f(\theta),v\right>  \\
 \, = \,
\left<\nabla f(\theta),v\right> - 
\left<\nabla f(\theta),\theta\right> \left<v,\theta\right>,
$$
let us extend it smoothly to all points $x \in G$ by
\be
\psi(x) = 
\left<\nabla f(x),v\right> - \left<\nabla f(x),x\right> \left<v,x\right>
\en
and write
\be
\nabla_S \psi(\theta) =
\nabla \psi(\theta) - \left<\nabla \psi(\theta),\theta\right> \theta.
\en

From (9.3) and (9.2) it follows that
\bee
\nabla \psi(x)
 & = &
\nabla \left<\nabla f(x),v\right> - 
\nabla \big(\left<\nabla f(x),x\right> \left<v,x\right>\big) \\
 & = &
f''(x)v - \left<\nabla f(x),x\right> v -
\nabla \big(\left<\nabla f(x),x\right>\big) \left<v,x\right>.
\ene
In addition, the function $u(x) = \left<\nabla f(x),x\right>$ has the 
gradient $\nabla u(x) = f''(x)x + \nabla f(x)$, so
$$
\nabla \psi(x) = f''(x)v - \left<\nabla f(x),x\right> v -
\big(f''(x)x + \nabla f(x)\big) \left<v,x\right>.
$$
Restricting this equality to the unit sphere and using for short 
the notation $P = P_{\theta^\perp}$, we get
\be
\nabla \psi(\theta) = f''(\theta) Pv - 
\left<\nabla f(\theta),\theta\right> v -
\left<v,\theta\right> \nabla f(\theta),
\en
which also implies
\be
\left<\nabla \psi(\theta),\theta\right>\theta = 
\left<f''(\theta) Pv,\theta\right>\theta - 
\left<\nabla f(\theta),\theta\right> \left<v,\theta\right>\theta -
\left<v,\theta\right> \left<\nabla f(\theta),\theta\right>\theta.
\en

Following (9.4), it remains to subtract (9.6) from (9.5). 
First note that
$$
f''(\theta) Pv - \left<f''(\theta) Pv,\theta\right>\theta = P f''(\theta) Pv
$$
which is deduced from the general formula
$Pw = w - \left<w,\theta\right>\theta$ with $w = f''(\theta) Pv$.
The equality $v - \left<v,\theta\right>\theta = Pv$ can be used for the 
second terms on the right of (9.5)-(9.6). Finally, for
the third terms we have
$$
\nabla f(\theta) - \left<\nabla f(\theta),\theta\right>\theta = 
P\,\nabla f(\theta) = \nabla_S f(\theta).
$$
Therefore, using the matrix $B$ from Proposition 7.1, the difference between 
(9.5) and (9.6) is exactly
\bee
P f''(\theta) Pv - \left<\nabla f(\theta),\theta\right> Pv  - 
\left<v,\theta\right> \nabla_S f(\theta)
 & = &
PBP - \left<v,\theta\right> \nabla_S f(\theta) \\
 & = &
f''_S(\theta) v - \left<v,\theta\right> \nabla_S f(\theta).
\ene
Thus,
\be
\nabla_S \left<\nabla_S f(\theta),v\right> = 
f''_S(\theta) v - \left<v,\theta\right> \nabla_S f(\theta),
\en
which is the desired equality (9.1).
\qed

\vskip10mm
\section{{\bf Appendix C. Second Order Modulus of Gradients}}
\setcounter{equation}{0}

\vskip2mm
\noindent
Let us give more details explaining Lemma 3.1. Recall that, by the very
definition of the second order modulus of the gradient,
\bee
|\nabla_S^2 f(\theta)|
 & = &
|\nabla_S \, |\nabla_S f(\theta)|\,| \\
 & = & 
\limsup_{\theta' \rightarrow \theta} 
\frac{|\,|\nabla_S f(\theta)|-|\nabla_S f(\theta')|\,|}{|\theta - \theta'|}, 
\qquad \theta \in S^{n-1}.
\ene

\vskip2mm
{\bf Proof of Lemma 3.1.}
First let us show that the function $|\nabla_S f|$ has a finite Lipschitz 
semi-norm. Since the first two spherical derivatives of $f$ are continuous 
and therefore bounded on the unit sphere, it follows from (9.7) that
$$
|\nabla_S \left<\nabla_S f(\theta),v\right>| \leq C, \qquad |v|=1,
$$
with some constant $C$ (independent of $\theta$ and $v$).
Hence, the function $\theta \rightarrow \left<\nabla_S f(\theta),v\right>$
has Lipschitz semi-norm at most $C$, so that
$$
|\left<\nabla_S f(\theta'),v\right> - \left<\nabla_S f(\theta),v\right>| \leq 
C \rho(\theta',\theta)
$$
for all $\theta, \theta' \in S^{n-1}$. Taking here the supremun over all 
unit vectors $v$ and applying the triangle inequality, we get
$$
\Big|\,|\nabla_S f(\theta')| - |\nabla_S f(\theta)|\,\Big| \leq
|\nabla_S f(\theta') - \nabla_S f(\theta)| \leq 
C \rho(\theta',\theta),
$$
which is the Lipschitz property (with constant $C$).

Next, to derive the required identity for the second order modulus of the 
gradient, we fix $\theta \in S^{n-1}$ and apply the identity (9.7) once 
more. By the definition of the spherical gradient, it yields the Taylor 
expansion up to the linear term,
\be
\left<\nabla_S f(\theta'),v\right> = \left<\nabla_S f(\theta),v\right> +
\left<V,\theta' - \theta\right> + o(|\theta' - \theta|)
\en
as $\theta' \rightarrow \theta$, $\theta' \in S^{n-1}$, where
$$
V = f''_S(\theta) v - \left<v,\theta\right> \nabla_S f(\theta).
$$
Moreover, by the Taylor formula in the integral form, and since any continuous
function on a compact metric space is uniformly continuous,
the $o$-term in (10.1) can be bounded by a quantity which is
independent of $v \in S^{n-1}$. That is,
$$
\sup_{v \in S^{n-1}} \,
|\left<\nabla_S f(\theta') - \nabla_S f(\theta),v\right> - 
\left<V,\theta' - \theta\right> \leq \ep(|\theta' - \theta|)
$$
with some function $\ep(t)$ such that $\ep(t) \rightarrow 0$ as
$t \rightarrow 0$. 

Now, let us rewrite (10.1) as
\be
\left<\nabla_S f(\theta'),v\right> = \left<\nabla_S f(\theta) + L,v\right>
 + o(|\theta' - \theta|),
\en
where
$$
\left<L,v\right> = \left<V,\theta' - \theta\right> =
\left<f''_S(\theta)v - 
\left<v,\theta\right> \nabla_S f(\theta),\theta' - \theta\right>,
$$
that is, with
\be
L = f''_S(\theta)(\theta' - \theta) - 
\left<\nabla_S f(\theta),\theta' - \theta\right>\theta.
\en
Taking an absolute value of both sides in (10.2) and turning to the 
supremum over all $v \in S^{n-1}$, we obtain that
\be
|\nabla_S f(\theta')| = |\nabla_S f(\theta) + L| + 
o(|\theta' - \theta|).
\en

Next, write
\be
|\nabla_S f(\theta) + L|^2 =
|\nabla_S f(\theta)|^2 + 2\left<\nabla_S f(\theta),L\right> + |L|^2.
\en
Since $\nabla_S f(\theta)$ is orthogonal to the vector $\theta$, 
we have from (10.3) that
$$
\left<\nabla_S f(\theta),L\right> = 
\left<\nabla_S f(\theta),f''_S(\theta)(\theta' - \theta)\right> =
\left<w,\theta' - \theta\right>,
$$
where 
$$
w = f''_S(\theta)\nabla_S f(\theta).
$$
Since also $|L|^2 = O(|\theta' - \theta|^2)$, (10.5) yields
$$
|\nabla_S f(\theta) + L|^2 =
|\nabla_S f(\theta)|^2 + 2\left<w,\theta' - \theta\right> + 
o(|\theta' - \theta|),
$$
and therefore in case $|\nabla_S f(\theta)|>0$,
$$
|\nabla_S f(\theta) + L| =
|\nabla_S f(\theta)| + |\nabla_S f(\theta)|^{-1}
\left<w,\theta' - \theta\right> + o(|\theta' - \theta|).
$$
Using this in (10.4), we find that
$$
|\nabla_S f(\theta')| - |\nabla_S f(\theta)| = 
|\nabla_S f(\theta)|^{-1}
\left<w,\theta' - \theta\right> + o(|\theta' - \theta|)
$$
and hence
\bee
\limsup_{\theta' \rightarrow \theta} 
\frac{|\,|\nabla_S f(\theta')|-|\nabla_S f(\theta)|\,|}{|\theta' - \theta|}
 & = &
|\nabla_S f(\theta)|^{-1}
\limsup_{\theta' \rightarrow \theta} 
\frac{|\left<w,\theta' - \theta\right>|}{|\theta' - \theta|} \\
 & = &
|\nabla_S f(\theta)|^{-1} P_{\theta^\perp} w.
\ene
Thus, by the definition, 
$$
|\nabla_S^2 f(\theta)| = |\nabla_S f(\theta)|^{-1} P_{\theta^\perp} w.
$$
But, as was noted before, the vector $w$ is always orthogonal to $\theta$. 
Therefore, $P_{\theta^\perp} w = w$, and we arrive at the required identity
$$
|\nabla_S^2 f(\theta)| = |\nabla_S f(\theta)|^{-1}\, w \quad
\big(|\nabla_S f(\theta)|>0\big).
$$

Finally, consider the remaining case $|\nabla_S f(\theta)| = 0$. Then
$L = f''_S(\theta)(\theta' - \theta)$, and (10.4) is simplified to
$$
|\nabla_S f(\theta')| = |L| + o(|\theta' - \theta|).
$$
Again, by the very definition, and using orthogonality of 
$f''_S(\theta)h$ to $\theta$,
\bee
|\nabla_S^2 f(\theta)|
 & = &
\limsup_{\theta' \rightarrow \theta} 
\frac{|\nabla_S f(\theta')|}{|\theta' - \theta|} \\
 & = &
\limsup_{\theta' \rightarrow \theta} 
\frac{|f''_S(\theta)(\theta' - \theta)|}{|\theta' - \theta|} \ = \
\limsup_{h \rightarrow 0, \ h \in \theta^\perp} 
\frac{|f''_S(\theta)h|}{|h|} \ = \ \|f_S''(\theta)\|.
\ene
\qed

\vskip10mm
\section{{\bf Appendix D. Laplacian}}
\setcounter{equation}{0}

\vskip2mm
\noindent
The Laplacian operator $\Delta_S f = {\rm Tr}\, f''_S$, acting in the
class of all $C^2$-smooth function $f$ on $S^{n-1}$, can be related
to the "spherical partial derivatives"
$
D_i f(\theta) = \left<\nabla_S f(\theta),e_i\right>,
$
where $e_1,\dots, e_n$ is the canonical basis in $\R^n$. Thus,
$$
\nabla_S f(\theta) = \sum_{i=1}^n D_i f(\theta)\, e_i.
$$
As the next partial derivatives, consider "second order" differential 
operators
$$
D_{ij} f = D_i (D_j f) = \left<\nabla_S \left<\nabla_S f,e_j\right>,e_i\right>, 
\qquad i,j = 1,\dots,n.
$$

\vskip5mm
{\bf Proposition 11.1.} {\it
$
\Delta_S = \sum_{i=1}^n D_{ii}.
$
}

\vskip5mm
In fact, any orthonormal basis in $\R^n$ could be used in place of $e_i$'s 
in the definition of $D_{ii}$, and the above statement will continue to hold.

\vskip5mm
{\bf Proof.} By (9.1), for all $v \in \R^n$,
$$
\nabla_S \left<\nabla_S f(\theta),v\right> = f''_S(\theta) v - 
\left<v,\theta\right> \nabla_S f(\theta).
$$
Hence
$$
D_{ii} f(\theta) = \left<f''_S(\theta) e_i,e_i\right> - 
\left<\theta,e_i\right> \left<\nabla_S f(\theta),e_i\right>
$$
and thus
$$
\sum_{i=1}^n D_{ii} f(\theta) \, = \,
{\rm Tr}\, f''_S(\theta) - \left<\nabla_S f(\theta),\theta\right>
 \, = \, {\rm Tr}\, f''_S(\theta).
$$
\qed

\vskip2mm
Let us now return to Lemma 4.5 with its identity
\be
\Delta_S f(\theta) =
\Delta f(\theta) - (n-1) \left<\nabla f(\theta),\theta\right> 
- \left<f''(\theta)\theta,\theta\right>.
\en
It can be obtained from the following explicit formula for the
derivatives $D_{ij}$. 

\vskip5mm
{\bf Lemma 11.2.} {\it If $f$ is $C^2$-smooth in 
an open neighborhood of $S^{n-1}$, then for all $\theta \in S^{n-1}$
and all $i,j=1,\dots,n$,
\begin{eqnarray}
D_{ij} f(\theta)
 & = &
\partial_{ij} f(\theta) - \theta_j\, \partial_i f(\theta) 
- \delta_{ij} \left<\nabla f(\theta),\theta\right> 
+ 2 \theta_i \theta_j \left<\nabla f(\theta),\theta\right> \nonumber \\
 & & \
- \theta_j \left<f''(\theta)\theta,e_i\right> - 
\theta_i \left<f''(\theta)\theta,e_j\right>
+ \theta_i \theta_j \left<f''(\theta)\theta,\theta\right>.
\end{eqnarray}
In particular,
\bee
D_{ii} f(\theta)
 & = &
\partial_{ii} f(\theta) - \left<\nabla f(\theta),\theta\right> -
\theta_i\, \partial_i f(\theta) 
+ 2 \theta_i^2 \left<\nabla f(\theta),\theta\right> \\
 & &
- \ 2\theta_i \left<f''(\theta)\theta,e_i\right> 
+ \theta_i^2 \left<f''(\theta)\theta,\theta\right>.
\ene
}

\vskip2mm
Summing the latter equality over all $i \leq n$, we arrive at (11.1).

Note that the operators $D_i$ and $D_j$ are not commutative, i.e. 
we do not have the identity $D_{ij} f = D_{ji} f$ in the entire
class $C^2$. Indeed, by (11.2),
$D_{ij} f(\theta) = D_{ji} f(\theta)$, if and only if
$\theta_i\, \partial_j f(\theta) = \theta_j\, \partial_i f(\theta)$.

\vskip5mm
{\bf Proof of Lemma 11.2.} Assume $f$ is $C^2$-smooth in the open region $G$.
Fix an index 
$j \leq n$ and consider the smooth function in $n$ real variables
\begin{eqnarray}
u(x) & = & 
\left<\nabla f(x),e_j\right> - 
\left<\nabla f(x),x\right>\left<x,e_j\right> \nonumber \\
 & = &
\partial_j f(x) - x_j \sum_{k=1}^n x_k\, \partial_k f(x), \qquad \qquad
x = (x_1,\dots,x_n) \in G.
\end{eqnarray}
In particular, $u(\theta) = D_j\, f(\theta)$ for $\theta \in S^{n-1}$ 
and therefore $D_{ij} f = D_i u$. Again following the definition of $D_i$, 
we have
\be
D_i\, u(x) = \partial_i u(x) - x_i \sum_{l=1}^n x_l\, \partial_l u(x).
\en
By (11.3),
\bee
\partial_i u(x)
 & = & 
\partial_{ij} f(x) - \delta_{ij} \sum_{k=1}^n x_k\, \partial_k f(x) - 
x_j\, \partial_i f(x) - x_j \sum_{k=1}^n x_k\, \partial_{ik} f(x), \\
\partial_l u(x)
 & = & 
\partial_{lj} f(x) - \delta_{lj} \sum_{k=1}^n x_k\, \partial_k f(x) - 
x_j\, \partial_l f(x) - x_j \sum_{k=1}^n x_k\, \partial_{lk} f(x).
\ene
Plugging these equalities in (11.4), we get
\bee
D_i\, u(x)
 & = &
\partial_{ij} f(x) - \delta_{ij} \sum_{k=1}^n x_k\, \partial_k f(x) - 
x_j\, \partial_i f(x) - x_j \sum_{k=1}^n x_k\, \partial_{ik} f(x) \\
 & & \
 - x_i \sum_{l=1}^n x_l\,
\Big[\partial_{lj} f(x) - \delta_{lj} \sum_{k=1}^n x_k\, \partial_k f(x) - 
x_j\, \partial_l f(x) - x_j \sum_{k=1}^n x_k\, \partial_{lk} f(x)\Big] \\
 & = &
\partial_{ij} f(x) - \delta_{ij} \sum_{k=1}^n x_k\, \partial_k f(x) - 
x_j\, \partial_i f(x) - x_j \sum_{k=1}^n x_k\, \partial_{ik} f(x) \\
 & & \
- x_i\, \sum_{l=1}^n x_l\,\partial_{lj} f(x) 
+ 2 x_i x_j\, \sum_{k=1}^n x_k\, \partial_k f(x)  
+ x_i x_j\, \sum_{l=1}^n \sum_{k=1}^n x_l x_k\, \partial_{lk} f(x). \\
\ene
In a bit more compact form,
\bee
D_i\, u(x)
 & = &
\partial_{ij} f(x) - x_j\, \partial_i f(x) 
- \delta_{ij} \left<\nabla f(x),x\right> 
+ 2 x_i x_j\, \left<\nabla f(x),x\right> \\
 & & \
- x_j \left<f''(x)x,e_i\right> - x_i\, \left<f''(x)x,e_j\right>
+ x_i x_j\, \left<f''(x)x,x\right>.
\ene
It remains to restrict this function to the sphere.
\qed

\vskip10mm
\section{{\bf Appendix E. Homogeneous functions}}
\setcounter{equation}{0}

\vskip2mm
\noindent
A function $F: \R^n \setminus \{0\} \rightarrow \R$ is called homogeneous
of order $d$ (where $d$ is a real number), or $d$-homogeneous, if it 
satisfies the relation
$$
F(\lambda x) = \lambda^d\, F(x), \qquad x \neq 0, \ 
\lambda > 0.
$$
Any such function is completely determined by its values on the unit sphere.
Alternatively, starting from a function $f$ on $S^{n-1}$, one may define
its unique $d$-homogeneous extension by putting
$$
F(x) = r^d f(r^{-1} x), \qquad r = |x|, \ \ x \neq 0.
$$

For example, if $f = 1$, then $F(x) = |x|^d$.

In this section, we collect several formulas for the derivatives of 
$d$-homogeneous functions. We will use the notations
$$
r = |x|, \qquad \theta = r^{-1} x = \frac{x}{|x|} \ \ \ \ (x \neq 0).
$$

\vskip5mm
{\bf Proposition 12.1.} {\it For the $d$-homogeneous extension 
$F(x) = r^d f(\theta)$ of a $C^1$-smooth function $f$ on $S^{n-1}$, 
we have that, for all $x \neq 0$,
\be
\nabla F(x) = r^{d-1} \big[d \, f(\theta)\theta + \nabla_S f(\theta)\big].
\en
}

\vskip2mm
This formula can be easily verified by the direct differentiation
(assuming that $f$ is defined and $C^1$-smooth in a neighborhood of the 
sphere), so we omit the proofs. 

For example, for the $1$-homogeneous extension $F(x) = r f(\theta)$, we have
$$
\nabla F(x) = f(\theta)\theta + \nabla_S f(\theta), \qquad
|\nabla F(x)|^2 = f(\theta)^2 + |\nabla_S f(\theta)|^2.
$$
In this particular case, such functions may be used, for example, to recover
the Poincar\'e inequality on the sphere on the basis of the Poincar\'e-type 
inequality for the Gaussian measure (which in turn has many elementary proof).

For the $0$-homogeneous extension $F(x) = f(\theta)$, we have
$$
\nabla F(x) = r^{-1}\,\nabla_S f(\theta), 
$$
and thus the usual (Euclidean) and spherical gradients coincide on the unit 
sphere: $\nabla F = \nabla_S f$ on $S^{n-1}$.

It is therefore interesting to know whether a similar identity
holds for the second derivative as well. The answer is negative,
although some relationship does exist.

\vskip5mm
{\bf Proposition 12.2.} {\it For the $d$-homogeneous extension 
$F(x) = r^d f(\theta)$ of a $C^2$-smooth function $f$ on $S^{n-1}$, 
we have for all $x \neq 0$ and $v \in \R^n$,
\begin{eqnarray}
F''(x)v
 & = &
r^{d-2} \Big[d(d-1) f(\theta) \left<v,\theta\right> \theta + 
d f(\theta) P_{\theta^\perp}v \nonumber \\
 & & + \ 
(d-1) \left<\nabla_S f(\theta),v\right> \theta +
(d-1) \left<v,\theta\right> \nabla_S f(\theta) + f''_S(\theta) v\Big].
\end{eqnarray}
}

\vskip2mm
In interesting particular cases $d=0,1$, (12.2) is
simplified. For the $0$-homogeneous extension, we have
$$
F''(x)v = r^{-2}\, \big[-\left<\nabla_S f(\theta),v\right> \theta -
\left<v,\theta\right> \nabla_S f(\theta) + f''_S(\theta) v\big],
$$
while for the $1$-homogeneous extension,
$$
F''(x) = r^{-1}\, \big[f(\theta) P_{\theta^\perp} + f''_S(\theta)\big].
$$

{\bf Proof.} From (12.1),
$$
\left<\nabla F(x),v\right> = r^{d-1} 
\big[d f(\theta) \left<v,\theta\right>  + \left<\nabla_S f(\theta),v\right>\big].
$$
We are in position to apply (12.1) once more, now with $d-1$ in place of $d$
and with 
$$
\psi(\theta) = 
d f(\theta) \left<v,\theta\right>  + \left<\nabla_S f(\theta),v\right>
$$ 
in place of $f$. It gives
\begin{eqnarray}
F''(x)v
 & = &
\nabla \left<\nabla F(x),v\right> \nonumber \\
 & = &
r^{d-2} \big[(d-1) \psi(\theta) \theta  + \nabla_S \psi(\theta)\big]
 \nonumber \\
 & = &
r^{d-2} \big[\,d(d-1) f(\theta) \left<v,\theta\right> \theta  + 
(d-1) \left<\nabla_S f(\theta),v\right> \theta +
\nabla_S \psi(\theta)\big].
\end{eqnarray}
To develop the last gradient, using
$\nabla_S \left<v,\theta\right> = Pv$ ($P = P_{\theta^\perp}$), 
first write
$$
\nabla_S \psi(\theta) = d \left<v,\theta\right> \nabla_S f(\theta) +
d f(\theta) Pv + \nabla_S \left<\nabla_S f(\theta),v\right>.
$$
In order to evaluate the last gradient, we apply the identity (9.1), 
which gives
$$
\nabla_S \psi(\theta) = (d-1) \left<v,\theta\right> \nabla_S f(\theta) +
d f(\theta) Pv + f''_S(\theta) v.
$$
Inserting this expression in (12.3), we arrive at the formula (12.2).
\qed

\vskip5mm
{\bf Corollary 12.3.} {\it For the $d$-homogeneous extension 
$F(x) = r^d f(\theta)$ of a $C^2$-smooth function $f$ on $S^{n-1}$, 
we have for all $x \neq 0$,
\be
\Delta F(x) \, = \,
r^{d-2} \big[\,d(n+d-2) f(\theta) + \Delta_S f(\theta)\big].
\en
}

\vskip2mm
In particular, $\Delta F(x) = r^{d-2} \Delta_S f(\theta)$
for the $0$-homogeneous extension $F(x) = f(\theta)$,
so the Euclidean and spherical Laplacians coincide on the unit sphere.
The same conclusion is also true when $d = 2-n$.

The identity (12.4) is well-known. It implies that, for any spherical 
harmonic $f$ on $S^{n-1}$ of degree $d$ (so that $\Delta F = 0$), we 
necessarily have $\Delta_S f = -d(n+d-2) f$, cf. e.g. [S-W].

\vskip5mm
{\bf Proof.} Applying (12.2) with $v = e_i$, we get
\bee
\left<F''(x)e_i.e_i\right>
 & = &
r^{d-2} \big[d(d-1) f(\theta) \left<\theta,e_i\right>^2 + 
d f(\theta) \left<P_{\theta^\perp}e_i,e_i\right>  \\
 & & + \ 
2(d-1) \left<\nabla_S f(\theta),e_i\right> \left<\theta,e_i\right> +
\left<f''_S(\theta)e_i,e_i\right>\big].
\ene
Here $P_{\theta^\perp}e_i = e_i - \left<\theta,e_i\right>\theta$, so
$$
\sum_{i=1}^n \left<P_{\theta^\perp}e_i,e_i\right> =
\sum_{i=1}^n \big(1 - \left<\theta,e_i\right>^2\big) = n-1.
$$
In addition,
$$
\sum_{i=1}^n \left<\nabla_S f(\theta),e_i\right> \left<\theta,e_i\right> =
\left<\nabla_S f(\theta),\theta\right> = 0.
$$
Hence,
$$
\Delta F(x) = 
\sum_{i=1}^n \left<F''(x)e_i.e_i\right> =
r^{d-2} \big[d(d-1) f(\theta) + 
d (n-1)\,f(\theta) + \Delta_S f''_S(\theta)\big].
$$
\qed

\vskip10mm
\section{{\bf Appendix F. Commutator of Laplacian and Gradient}}
\setcounter{equation}{0}

\vskip2mm
\noindent
If a function $f$ is defined and $C^3$-smooth in an open region of $\R^n$,
then, for any $v \in \R^n$,
\be
\Delta \left<\nabla f(x),v\right> = \left<\nabla \Delta f(x),v\right>
\en
throughout the region. In a more compact form,
$\Delta \nabla = \nabla \Delta$, that is, these two operators -- the
Euclidean Laplacian and the Euclidean gradient -- commute.
However, due to curvature of $S^{n-1}$, this is no longer true for 
the spherical Laplacian and the spherical gradient which may be seen from 
the formula for the commutator given in Lemma 4.3.
In a more compact vector form, this formula may be written as
$$
\Delta_S \nabla_S f(\theta) - \nabla_S \Delta_S f(\theta) \, = \,  
(n-3)\,\nabla_S f(\theta) - 2\Delta_S f(\theta)\,\theta.
$$
This identity may also be rewritten component-wise in terms of the operators 
$D_i$ as
$$
\Delta_S D_i f(\theta) - D_i\, \Delta_S f(\theta) \, = \,  
(n-3) D_i f(\theta) - 2\left<\theta,e_i\right>\Delta_S f(\theta).
$$

\vskip5mm
{\bf Proof of Lemma 4.3.} By Proposition 12.1 and Corollary 12.3, 
for any $C^3$-smooth $d$-homogeneous function $u$ on $\R^n \setminus \{0\}$, 
for all $x \neq 0$,
\begin{eqnarray}
\nabla u(x) 
 & = &
r^{d-1} \big[\,d \, u(\theta)\theta + \nabla_S u(\theta)\big], \\ 
\Delta u(x) 
 & = &
r^{d-2} \big[\,d(n+d-2)\, u(\theta) + \Delta_S u(\theta)\big],
\end{eqnarray}
where $r = |x|$ and $\theta = r^{-1} x$.

The identity (13.1) will be used with the $0$-homogeneous extension 
$F(x) = f(\theta)$, $x \neq 0$. Being restricted to the points lying 
on the unit sphere, it becomes
\be
\Delta \left<\nabla F(\theta),v\right> = 
\left<\nabla \Delta F(\theta),v\right>.
\en
In that case, $\nabla_S f(\theta) = \nabla F(\theta)$, so
$$
\Delta_S \left<\nabla_S f(\theta),v\right> =
\Delta_S \left<\nabla F(\theta),v\right>.
$$
Moreover, the function $u(x) = \left<\nabla F(x),v\right>$ is 
$(-1)$-homogeneous, and we may apply (13.3) with $d=-1$. Again,
being restricted to the unit sphere, this identity becomes
$$
\Delta u(\theta) = -(n-3)\, u(\theta) + \Delta_S u(\theta),
$$
so
\begin{eqnarray}
\Delta_S \left<\nabla_S f(\theta),v\right>
 \ = \ \Delta_S u(\theta)
 & = & \Delta u(\theta) + (n-3)\, u(\theta) \nonumber \\
 & = &
\Delta \left<\nabla F(\theta),v\right> +  (n-3)\,\left<\nabla F(\theta),v\right>.
\end{eqnarray}

On the other hand, the function $u(x) = \Delta F(x)$ is 
$(-2)$-homogeneous, and we may apply (13.2) with $d=-2$. It gives
$$
\nabla u(\theta) = -2 u(\theta)\,\theta + \nabla_S u(\theta).
$$
Since $\Delta F$ coincides with $\Delta_S f$ on $S^{n-1}$, we get that
\begin{eqnarray}
\left<\nabla_S \Delta_S f(\theta),v\right>
 & = & 
\left<\nabla_S\, u(\theta),v\right>  \nonumber \\
 & = & 
\left<\nabla u(\theta),v\right> + 2 u(\theta)\left<\theta,v\right> \nonumber \\
 & = &
\left<\nabla \Delta F(\theta),v\right> + 2 \Delta F(\theta)\left<\theta,v\right>.
\end{eqnarray}

It remains to subtract (13.6) from (13.5) and apply (13.4), which 
leads to
$$
\Delta_S \left<\nabla_S f(\theta),v\right> - 
\left<\nabla_S \Delta_S f(\theta),v\right> \, = \, 
(n-3)\left<\nabla F(\theta),v\right> - 2\left<\theta,v\right>\Delta F(\theta).
$$
But $\nabla F(\theta) = \nabla_S f(\theta)$ and
$\Delta F(\theta) = \Delta_S f(\theta)$.
\qed

\vskip10mm
\section{{\bf Appendix G. Integrals Involving Laplacian}}
\setcounter{equation}{0}

\vskip2mm
\noindent
Many integrals involving the spherical Laplacian can be evaluated with the help of 
the classical formula
\be
\int \left<\nabla_S f,\nabla_S g\right> d\sigma_{n-1} = -
\int f \Delta_S g\,d\sigma_{n-1},
\en
which actually may be taken as an equivalent definition of the operator $\Delta_S$.
It yields the following characterization of the orthogonality to linear
functions in terms of the spherical gradient.

\vskip5mm
{\bf Proposition 14.1.} {\it For any smooth function $f$ on $S^{n-1}$,
\be
\int f(\theta) \theta\,d\sigma_{n-1}(\theta) = 
\frac{1}{n-1} \int \nabla_S f(\theta)\,d\sigma_{n-1}(\theta).
\en
In particular, $f$ is orthogonal to all linear functions in 
$L^2(S^{n-1})$, if and only if all linear forms 
$\left<\nabla_S f(\theta),v\right>$ have $\sigma_{n-1}$-mean zero.
}

\vskip5mm
{\bf Proof.} The linear function  
$g(\theta) = \left<v,\theta\right>$ has the spherical gradient and
respectively the spherical Laplacian
$$
\nabla_S g(\theta) = P_{\theta^\perp} v, \qquad
\Delta_S g(\theta) = - (n-1) \left<v,\theta\right>.
$$
In this case, (14.1) becomes
$$
\int \left<\nabla_S f(\theta),P_{\theta^\perp} v\right> d\sigma_{n-1}(\theta) = 
(n-1) \int f(\theta) \left<v,\theta\right> d\sigma_{n-1}(\theta).
$$
But $\nabla_S f(\theta)$ is orthogonal to $\theta$, so, 
$
\left<\nabla_S f(\theta),P_{\theta^\perp} v\right> =
\left<\nabla_S f(\theta),v\right>.
$
It follows that
$$
\int \left<\nabla_S f,v\right> d\sigma_{n-1} = 
(n-1) \int f(\theta) \left<v,\theta\right> d\sigma_{n-1}(\theta)
$$
which is the required identity (14.2). 
\qed

\vskip5mm
In fact, the formula (14.1) and the one of Lemma 4.4 may be extended to a more 
general Green-type formula with weights.

\vskip5mm
{\bf Proposition 14.2.} {\it For all $C^2$-smooth functions $f,g$ and any
$C^1$-smooth function $w$ on $S^{n-1}$,
\be
\int \left<\nabla_S f,\nabla_S g\right>\,w\,d\sigma_{n-1} = -
\int f\Delta_S g\,w\,d\sigma_{n-1} -
\int f\left<\nabla_S g,\nabla_S w\right> d\sigma_{n-1}.
\en
}

\vskip2mm
When $w=1$, we return to (14.1). In the case of the linear weight 
$w(\theta) = \left<v,\theta\right>$, let us recall that 
$\nabla_S w(\theta) = P_{\theta^\perp} v = v - \left<v,\theta\right>\theta$ 
and that $\nabla_S g(\theta)$ is orthogonal to $\theta$. Hence, (14.3)
is simplified to
\bee
\int \left<\nabla_S f(\theta),\nabla_S g(\theta)\right>\left<v,\theta\right> 
d\sigma_{n-1}(\theta)
 & = &
-
\int f(\theta)\Delta_S g(\theta)\left<v,\theta\right> d\sigma_{n-1}(\theta) \\
 & & - \ 
\int f(\theta)\left<\nabla_S g(\theta),v\right>\,d\sigma_{n-1}(\theta),
\ene
which is the statement of Lemma 4.4.

\vskip5mm
{\bf Proof of Proposition 14.2.}
Using the canonical basis in $\R^n$, first write
\begin{eqnarray}
\int \left<\nabla_S f,\nabla_S g\right> w\,d\sigma_{n-1}
 & = &
\sum_{i=1}^n \int \left<\nabla_S f,e_i\right> 
\left<\nabla_S g,e_i\right> w\,d\sigma_{n-1} \nonumber \\
 & = &
\sum_{i=1}^n \left<\int
\left<\nabla_S g,e_i\right> w\, \nabla_S f\,d\sigma_{n-1},e_i\right>.
\end{eqnarray}
Applying the general identity
$\nabla_S (\varphi \psi) =  \psi\nabla_S \varphi + \varphi\nabla_S \psi$, 
let us represent the last vector integral in (14.4) as $I_i + J_i$, where
$$
I_i = - \int f\, \nabla_S \big(\left<\nabla_S g,e_i\right> w\big)\,d\sigma_{n-1},
\qquad
J_i = \int \nabla_S \big(f \left<\nabla_S g,e_i\right> w\big)\,d\sigma_{n-1}.
$$
Again by the same identity,
\bee
I_i
 & = &
- \int f\, \Big(w\,\nabla_S \left<\nabla_S g,e_i\right>
+ \left<\nabla_S g,e_i\right> \nabla_S w\Big)\,d\sigma_{n-1} \\
 & = &
- \int f\, w\, D_{ii} g\,d\sigma_{n-1} \ e_i
- \int f \left<\nabla_S g,e_i\right> \nabla_S w\,d\sigma_{n-1}.
\ene
Hence,
\bee
\sum_{i=1}^n \left<I_i,e_i\right> 
 & = & 
- \int f w\, \Delta_S g\,d\sigma_{n-1} - \sum_{i=1}^n 
\int f \left<\nabla_S g,e_i\right> \left<\nabla_S w,e_i\right>\,d\sigma_{n-1} \\
 & = &
- \int fw\, \Delta_S g\,d\sigma_{n-1} -
\int f\left<\nabla_S g,\nabla_S w\right>\,d\sigma_{n-1},
\ene
which is exactly the desired expression appearing on the right-hand side of (14.3).

On the other hand, by Proposition 14.1,
$$
J_i = (n-1) 
\int f(\theta) \left<\nabla_S g(\theta),e_i\right> w(\theta) \theta\,d\sigma_{n-1}(\theta),
$$
implying that
\bee
\sum_{i=1}^n \left<J_i,e_i\right> 
 & = &
(n-1) \sum_{i=1}^n
\int f(\theta) \left<\nabla_S g(\theta),e_i\right> 
w(\theta) \left<\theta,e_i\right>d\sigma_{n-1}(\theta) \\
 & = &
(n-1) 
\int f(\theta) w(\theta) \left<\nabla_S g(\theta),\theta\right> d\sigma_{n-1}(\theta) \ = \ 0.
\ene
Thus, only the integrals $I_i$'s contribute in the sum (14.4).
\qed

\end{document}